\PassOptionsToPackage{usenames,dvipsnames}{xcolor}
\documentclass[a4paper, twocolumn=off, twoside=on]{scrartcl} 

\usepackage[hidelinks]{hyperref}
\usepackage{geometry}
\usepackage{orcidlink}

\PassOptionsToPackage{usenames,dvipsnames}{xcolor}
\usepackage[english]{babel}
\usepackage{amsmath}
\usepackage{amsthm}
\usepackage{amsfonts}
\usepackage{amssymb}
\usepackage{siunitx}
\usepackage{float}
\usepackage{caption}
\usepackage{subcaption}
\usepackage{xfrac}
\usepackage{booktabs}
\usepackage{xcolor}
\usepackage[only,llbracket,rrbracket]{stmaryrd}
\usepackage{extdash}
\usepackage{iftex}

\usepackage{enumitem}
\setlist{noitemsep, left=1em, rightmargin=1em}
\setitemize[2]{label=\textperiodcentered, rightmargin=1em}

\usepackage{siunitx}
\usepackage[fixamsmath]{mathtools}
\mathtoolsset{showonlyrefs=true}

\theoremstyle{definition}
\newtheorem{definition}{Definition}

\theoremstyle{plain}

\theoremstyle{remark}
\newtheorem{remark}[definition]{Remark}

\addto\extrasenglish{}

\usepackage[proportional]{libertinus}
\usepackage{libertinust1math}
\usepackage[cal=boondox,scr=boondox, bb=ams]{mathalpha}

\newcommand{\symbf}[1]{\mathbfit{#1}}
\newcommand{\symsf}[1]{\mathsf{#1}}
\newcommand{\symbfsf}[1]{\mathsfbf{#1}}
\newcommand{\symbfup}[1]{\mathbf{#1}}
\newcommand{\symbb}[1]{\mathbb{#1}}
\newcommand{\symbfcal}[1]{\mathbfcal{#1}}
\newcommand{\symcal}[1]{\mathcal{#1}}

\KOMAoptions{parskip=half}
\KOMAoptions{abstract=true}

\newenvironment{frontmatter}{}{}

\addtokomafont{section}{}
\addtokomafont{titlehead}{}
\addtokomafont{subject}{\normalfont\sffamily\normalsize}
\addtokomafont{title}{\normalfont\sffamily\normalfont\sffamily\Large}
\addtokomafont{subtitle}{\normalfont\sffamily\large}
\addtokomafont{author}{\normalfont\sffamily\large}
\addtokomafont{date}{\normalfont\sffamily\large}
\addtokomafont{publishers}{\normalfont\sffamily\normalsize}
\setkomafont{paragraph}{\normalfont\sffamily\normalsize\itshape}

\usepackage{authblk}

\usepackage[en-US]{datetime2}
\DTMlangsetup{showdayofmonth=false}

\setenumerate[1]{label={\arabic*}.,  leftmargin=*}
\setitemize{leftmargin=*}
\setlist{rightmargin=1em}

\usepackage[clines, plainfootsepline, automark]{scrlayer-scrpage}
\clearpairofpagestyles

\ihead[]{}
\cehead{J. Magiera, C. Rohde}
\cohead{\csname @title\endcsname}
\ohead[]{\pagemark}
\ifoot[]{}
\ofoot[]{}

\setkomafont{pagehead}{\normalfont\footnotesize\itshape}
\setkomafont{pagefoot}{\normalfont\footnotesize}
\setkomafont{pagenumber}{\normalfont\normalsize\itshape}

\usepackage{xifthen}
\usepackage{xparse}

\ifthenelse{\isempty{#1}}{}{}

\usepackage[many]{tcolorbox}
\tcbuselibrary{theorems}
\tcbset{
breakable, 
enhanced jigsaw, 
lines before break=20,
colback=white,
colframe=black!15, 
coltitle=black, 
colbacktitle=black!15, 
fonttitle=\bfseries\sffamily, arc=0mm,
subtitle style={colback=black!15,fontupper={\sffamily\color{black}}}
}
\tcbset{shield externalize}

\newtcolorbox[use counter*=definition]{boxalgorithm}[2][]{
list text={#2},
  title={Algorithm \thetcbcounter: #2},
  #1
}

\newcommand{\pliq}{-}
\newcommand{\pvap}{+}
\newcommand{\phasel}{\pliq}
\newcommand{\phaser}{\pvap}

\newcommand{\phaserl}{\pm}
 \newcommand{\datx}{\symbfsf{x}} \newcommand{\daty}{\symbfsf{y}} \newcommand{\datz}{\symbfsf{z}} \newcommand{\cres}{CRes}  \newcommand{\ipmm}{IPMM}
\newcommand{\mdshort}{MD} 
\newcommand{\eosshort}{EOS}
\newcommand{\eosplshort}{EOS}

\newcommand{\nparticles}{\ensuremath{N_{\mathrm{p}}}} 

\newcommand{\mdtimevar}{\tau} 
\newcommand{\mdinterface}{\Gamma} 
\newcommand{\mddomain}{\Omega} 
\newcommand{\mdtimeend}{\mdtimevar_{\mathrm{end}}} 
\newcommand{\mdstepend}{N_{\mathrm{end}}} 

\newcommand{\continterface}{\symbfup{\Gamma}} 
\newcommand{\contdomain}{\symbfup{\Omega}} 

\newcommand{\nnparams}{\symbfup{\theta}} 
\newcommand{\datasetvar}{D} 

\newcommand{\frequencyvar}{\chi}

\newcommand{\networkF}{\symsf{R}} \newcommand{\networkG}{\symsf{G}} 
\newcommand{\ipfv}{FV-IPMM}

\newcommand{\pluseq}{\mathrel{+}=}
\newcommand{\minuseq}{\mathrel{-}=}

\newcommand{\setsep}{\allowbreak:\allowbreak}

\newcommand{\rmicro}{R}  
\newcommand{\rmd}{\ensuremath{R_{\mathrm{MD}}}} \newcommand{\genericsolversym}{\symcal{R}}

\newcommand{\PP}{\mathcal{P}}
\newcommand{\pp}{P}

\newcommand{\uu}{\vec{u}} 
\newcommand{\UU}{\vec{U}} 

\renewcommand{\vec}[1]{\symbf{#1}}

\newcommand{\mat}[1]{\symbfsf{#1}}

\newcommand{\bI}{\symbb{I}}

\newcommand{\bN}{\symbb{N}}

\newcommand{\bR}{\symbb{R}}
\newcommand{\bS}{\symbb{S}}

\newcommand{\cC}{\mathcal{C}}

\newcommand{\cK}{\mathcal{K}}

\newcommand{\cR}{\mathcal{R}}
\newcommand{\cS}{\mathcal{S}}
\newcommand{\cT}{\mathcal{T}}
\newcommand{\cU}{\mathcal{U}}
\newcommand{\cV}{\mathcal{V}}

\newcommand{\bcK}{\symbfcal{K}}

\newcommand{\va}{\vec{a}}

\newcommand{\vf}{\vec{f}}

\newcommand{\vm}{\vec{m}}
\newcommand{\vn}{\vec{n}}

\newcommand{\vp}{\vec{p}}

\newcommand{\vr}{\vec{r}}
\newcommand{\vs}{\vec{s}}

\providecommand{\vv}{}
\renewcommand{\vv}{\vec{v}}

\newcommand{\vx}{\vec{x}}
\newcommand{\vy}{\vec{y}}

\newcommand{\vF}{\vec{F}}

\newcommand{\vV}{\vec{V}}

\DeclareMathOperator{\dd}{d}
\DeclareMathOperator{\ddd}{d\!}
\DeclareMathOperator{\DD}{D}

\DeclareMathOperator*{\argmax}{arg\,max}

\DeclarePairedDelimiter\abs{\lvert}{\rvert}

\DeclarePairedDelimiter\norm{\lVert}{\rVert}

\DeclarePairedDelimiter{\jump}{\llbracket}{\rrbracket}
\DeclarePairedDelimiterX{\mean}[1]{\{}{\}}{\mkern-3.0mu\delimsize\{ {#1} \delimsize\}\mkern-3.0mu}

\usepackage{csquotes}
\usepackage[backend=biber, style=numeric, sortcites=true, citestyle=numeric-comp, giveninits=true, sorting=nyt, doi=true, isbn=false, url=false, maxcitenames=3, maxbibnames=10, minalphanames=3, block=none, date=year]{biblatex}
\AtEveryBibitem{\clearlist{language}\clearfield{pages}\clearfield{pagetotal}\clearfield{file}\clearlist{address}\clearlist{location}\clearlist{venue}}
    
\renewbibmacro{in:}{}

\begin{document}

\begin{frontmatter}

\title{A Molecular--Continuum Multiscale Model for Inviscid Liquid--Vapor Flow with Sharp Interfaces}

\author{Jim Magiera\,\orcidlink{0000-0001-9807-0784}}
\author{Christian Rohde\,\orcidlink{0000-0001-9183-5094}}
\affil{University of Stuttgart, Institute of Applied Analysis and Numerical Simulation}
\date{}
\maketitle

\begin{abstract}
\noindent{}The dynamics of compressible liquid--vapor flow depends sensitively on the microscale behavior at the phase boundary. 
  We consider a sharp-interface approach, and propose a multiscale model to describe liquid--vapor flow accurately, without imposing ad-hoc closure relations on the continuum scale.
  The multiscale model combines the Euler equations on the continuum scale with molecular-scale particle simulations that govern the interface motion. 
  We rely on an interface-preserving moving mesh finite volume method 
  to discretize the continuum-scale sharp-interface flow in a conservative manner.  
  Computational efficiency, while preserving physical properties, is achieved by a surrogate solver for the interface dynamics based on constraint-aware neural networks.
  \newline 
  The multiscale model is presented in its general form, and applied to regimes of temperature-dependent liquid--vapor flow which have not been accessible before. 
\end{abstract}
\end{frontmatter}

\section{Introduction}

In two-phase flow the fluid properties vary strongly among their different phases. 
This makes modeling two-phase flow quite challenging, as for example on a numerical level the mass transfer for large density jumps has to be resolved accurately. 
Depending on the spatial scale, there are several approaches to model two-phase flow.
Our focus lies on a spatial scale where all phase boundaries are explicitly represented.
In this work, we consider a sharp-interface approach for ideal fluid flow governed by the Euler equations as a system of first-order conservation laws. 
Thus, the fluid domain \(\contdomain \subset \bR^d\), \(d \in \bN\), is split into two subdomains: the liquid-phase domain \(\contdomain_{\pliq}\) and the vapor-phase domain \(\contdomain_{\pvap}\), 
and both domains are connected via a codimension 1 phase boundary \(\continterface(t)\).
In each of these domains only the respective phase is present, and both phases are coupled via the free interface using the Rankine--Hugoniot conditions. 
To determine the interface motion, further conditions have to be prescribed at the interface. 
Classically this is done by imposing an algebraic relation to fix the entropy dissipation rate \cite{truskinovsky:kinks:1993}. 
However, determining the explicit form of this relation remains a largely unsolved issue for complex flow regimes.
We model the interface flow dynamics on a more fundamental level and thus avoid prescribing some ad-hoc closure relations at the interface.
Precisely, we simulate two-phase flow directly on the molecular scale using molecular dynamics (\mdshort{}). 
This is attractive, because each molecule has the same properties and parameters regardless of the fluid phase it is part of.
With this ansatz we develop a multiscale model that describes interfacial mass, momentum, and energy transfer, as well as the motion of the phase boundary by an \mdshort{} interface solver. 

We note that there is no rigorous mathematical theory that allows us to directly and analytically bridge the scales from \mdshort{} systems with pairwise interactions to continuum-scale Euler equations.
Yet, it is still possible to use classical, statistical averaging to infer continuum-scale quantities from microscale particle simulations. 
In that way, we transfer the microscale behavior at fluid interfaces to the continuum scale via a data-based approach. 
This approach can be categorized as a heterogeneous multiscale method (HMM)
to couple atomistic- and continuum-scale models \cite{e.engquist.ea:heterogeneous:2007}. 
Moreover, we present a
finite volume method
on interface-preserving moving meshes (\ipfv{}-method)
which enables us to include the molecular-scale dynamics directly at the phase boundary, while resolving the sharp interface directly within the mesh.
The corresponding moving mesh-algorithms are presented in \cite{alkaemper.magiera.ea:interface:2021}.
\newline
Finally, to overcome the computational costs of expensive \mdshort{} simulations, we apply a surrogate for the \mdshort{} interface solver. 
This surrogate is based on constraint-aware neural networks \cite{magiera.ray.ea:constraint:2020}, to achieve computational efficiency, while preserving the mass-conservation of the numerical discretization. 

Two-phase fluid flow models with sharp interfaces have already been widely investigated, see e.g. 
\cite{saurel.petitpas.ea:modelling:2008,faccanoni.kokh.ea:modelling:2012,ma.lv.ea:entropy:2017,ghazi.james.ea:nonisothermal:2021} 
for various modeling approaches. 
Common discretization techniques involve for example cut-cell approaches
\cite{pan.han.ea:conservative:2018,shen.ren.ea:3d:2020} and 
interface tracking using level-sets \cite{fechter.munz.ea:sharp:2017,long.cai.ea:accelerated:2021}.
Compared to most of these, our discretization method does not require interface reconstruction steps, nor an extra scalar field for phase-tracking. 
\newline 
\mdshort{} simulations are a widely used tool in phase-transition applications.
However, the integration into continuum-scale dynamics for compressible liquid--vapor flow has been rarely done.
A direct comparison of \mdshort{} simulations and continuum-scale liquid--vapor flow has been carried out in e.g. \cite{frezzotti.barbante:simulation:2020,hitz.joens.ea:comparison:2021}.
In contrast to our contribution, the microscopic flow dynamics have not been directly integrated in the continuum-scale simulations over a whole range of flow regimes. 
\newline 
In recent years, a lot of attention was on the application of neural networks in the context of hyperbolic conservation laws. 
For example in \cite{patel.manickam.ea:thermodynamically:2022} physics-informed neural networks were designed for conservation laws. 
In \cite{bezgin.schmidt.ea:data:2021} neural networks are applied to infer reconstruction weights and dissipation terms for non-classical shock-waves. 
In our contribution, we focus on the application of neural networks as a surrogate solver, that assist in bridging the scales from the molecular to the continuum scale, while retaining important physical constraints such as mass conservation.

\paragraph*{Outline}
First, we introduce in Section~\ref{sec:generic_cl} a general sharp-interface two-phase model which is based on conservation laws. 
At the end of this section, we introduce a compressible liquid--vapor flow model, that serves as an explicit instance of the general model.
In Section~\ref{sec:interface_solver}, we introduce the concept of a microscale interface solver that is used to describe the phase boundary dynamics. 
Thereafter, we present an interface solver that is based on \mdshort{} simulations. 
For actual simulations, we do not employ the microscale interface solver directly, but instead use a surrogate solver that is described in Section~\ref{sec:surrogate}.
In Section~\ref{sec:ipvm_algorithm}, we provide the continuum-scale discretization scheme.
It uses a moving-mesh method \cite{alkaemper.magiera.ea:interface:2021} to resolve the sharp-interface directly within the mesh.  
In Section~\ref{sec:multiscale}, we present the multiscale model in a general form and provide a specific instance for the temperature-dependent two-phase flow model. 
In Section~\ref{sec:results:noniso}, we show numerical simulations results for this model. 
Finally, we discuss our findings and conclusions in Section~\ref{sec:conclusion}.

The content of this contribution is part of the PhD thesis \cite{magiera:molecular:2021}, preliminary parts are published in \cite{magiera.rohde:analysis:2021}.

 \section{Continuum-Scale Sharp-Interface Model for Two-Phase Flow}
\label{sec:generic_cl}

We begin by introducing a general sharp-interface model for two-phase phase flow, which is based on first-order conservation laws. 

\subsection{The General Two-Phase Flow Model}
\label{sec:generic_two_phase}

\begin{figure}[tb]
\centering
\includegraphics[width=0.2\columnwidth]{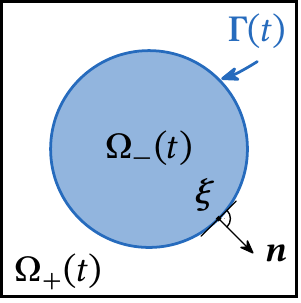}
\caption{
Sketch of the domain \(\contdomain\), that is divided by the phase boundary \(\continterface(t)\) into a \(\phasel\)-phase domain \(\contdomain_{\phasel}(t)\) and a surrounding \(\phaser\)-phase domain \(\contdomain_{\phaser}(t)\) at time \(t\). 
}
\label{fig:droplet_domain_iso}
\end{figure}

\paragraph*{Geometric Setting}
The central assumption for the sharp-interface approach is that the fluid domain \(\contdomain \subset \bR^{d}\), \(d \in \bN\), is divisible into two disjoint, open and \(d\)-dimensional subdomains \(\contdomain_\phasel(t)\), \(\contdomain_\phaser(t)\), fulfilling \((\overline{\contdomain_\phasel(t)} \cup \overline{\contdomain_\phaser(t)})^{\circ} = \contdomain\) for each \(t \geq 0\). 
Here, each bulk domain \(\contdomain_\phaserl(t)\) is filled with one respective fluid phase only, and the dividing hypersurface \(\continterface(t) = \overline{\contdomain_\phasel(t)} \cap \overline{\contdomain_\phaser(t)}\) is called the 
(sharp) interface, see Figure~\ref{fig:droplet_domain_iso} for a sketch of a typical setting.
For a point \(\vec{\xi} \in \continterface(t)\) the normal vector pointing in direction of the \(\phaser\)-phase is denoted by \(\vn = \vn(\vec{\xi},t) \in \bS^{d-1}\). 
The speed of the phase boundary \(\continterface(t)\) in normal direction is denoted by \(s = s(\vec{\xi},t) \in \bR\). 

\paragraph*{Bulk Equations}
We assume that in both bulk phases \(\contdomain_{\phaserl}(t)\)  the same set of first-order conservation laws are fulfilled, i.e.
\begin{align} 
\label{eq:generic_cl_2p}
\begin{aligned}
 \partial_t \UU^{\phaserl} + \nabla \cdot \vf(\UU^{\phaserl}) &= \vec{0},
 && \text{ in } \contdomain_\phaserl(t) \text{ for } t \in (0,t_{\mathrm{end}}), \\
 \UU^{\phaserl}(\:\cdot\:,0) &= \UU^{\phaserl}_{0} 
 && \text{ in } \contdomain_\phaserl(0),
 \end{aligned}
\end{align} 
with \(\UU^{\phaserl} \colon \contdomain_{\phaserl}(t) \times [0,t_{\mathrm{end}}) \to \cU\) being the unknown fluid states, with the open set \(\cU \subset \bR^m\), \(m \in \bN\), denoting the state space of the system. 
The flux function is given by \(\vf \in C^1(\cU, \bR^{m \times d})\), and the divergence operator in \eqref{eq:generic_cl_2p} is understood to act on each row of \(\vf(\UU^{\phaserl})\).
The initial conditions are given by \(\UU^{\phaserl}_{0}  \colon \contdomain(t) \to \cU\).
For a more compact notation, we write
\begin{align}
  \label{eq:compact_phase_notation}
  \UU(\vx, t) \coloneqq 
  \begin{cases}
    \UU^\phaser(\vx,t), & \text{ if } \vx \in \contdomain_\phaser(t), t \in [0,t_{\mathrm{end}}),  \\
    \UU^\phasel(\vx,t), & \text{ if } \vx \in \contdomain_\phasel(t), t \in [0,t_{\mathrm{end}}).
  \end{cases}  
\end{align}

We define the flux in direction of an arbitrary \(\vn \in \bS^{d-1}\) by 
\begin{align}
  \label{eq:generic_directional_flux}
  \vF(\UU ; \vn) \coloneqq \vf(\UU) \cdot \vn.
\end{align}
The system of conservation laws \eqref{eq:generic_cl_2p} is called hyperbolic, if for all \(\UU \in \cU\) and \(\vn \in \bS^{d-1}\) there exist \(m\) real eigenvalues 
\(\lambda_1(\UU ; \vn)\),\ldots, \(\lambda_m(\UU ; \vn)\), 
and \(m\) corresponding linear independent eigenvectors 
\(\vr_1(\UU ; \vn)\),\ldots, 
\(\vr_m(\UU ; \vn)\) of the flux Jacobian \(\DD\!\vF(\UU ; \vn) \colon \allowbreak \cU \to  \bR^{m\times m}\).

\paragraph*{Interface Conditions}
Let \(\vn = \vn(\vec{\xi}, t)\) be  the outer unit normal at \(\vec{\xi} \in \continterface(t)\).
We define the jump operator \(\jump{\cdot}\) at the interface \(\continterface\) 
for a sufficiently smooth function \(a \colon \bR^d \to \bR\) by
\begin{align}
  \label{eq:jump_mean}
  \jump{a} &\coloneqq \overline{a}_\phaser - \overline{a}_\phasel, 
  \text{ with }
  \overline{a}_\phaserl \coloneqq \lim_{\substack{\varepsilon \to 0 \\ \varepsilon > 0}}  a(\vec{\xi} \pm \varepsilon \vn).
\end{align}
To ensure conservation across discontinuities the Rankine--Hugoniot conditions 
\begin{align} 
  \label{eq:generic_rh}
s\,\jump{\UU} = \jump{\vF(\UU;\vn)}
\end{align}
must be fulfilled \cite{dafermos:hyperbolic:2016} .
Typically, for phase transition problems, further algebraic jump conditions are formulated. 
Those are classically formed by kinetic relations \cite{truskinovsky:kinks:1993} of the form 
\begin{align}
  \label{eq:generic_kinetic_relation}
  \bcK(\UU_{\phasel}, \UU_{\phaser}) = 0,
\end{align}
for a prescribed function \(\bcK \colon \PP_{\phasel} \times \PP_{\phaser} \to \bR\) that is usually associated with driving forces at the interface. 
Together with appropriate boundary conditions, the system \eqref{eq:generic_cl_2p}, with the jump conditions \eqref{eq:generic_rh} and  \eqref{eq:generic_kinetic_relation}, is a free boundary value problem for \(\UU^{\phaserl}\) and the interface \(\continterface\).

\paragraph*{Phase States}
Mathematically, the different phases are expressed as connected subsets in the state space \(\cU\) for which the system \eqref{eq:generic_cl_2p} is hyperbolic.
This, of course, depends heavily on the flux function \(\vf\) and the state space \(\cU\) itself.
For most real world fluids (as well as two-phase model fluids) the state space \(\cU\) contains (at least) two disjoint domains \(\PP_\phasel, \PP_\phaser \subset \cU\), in such a way that for all \(\UU \in \PP_\phasel \cup \PP_\phaser\) the system \eqref{eq:generic_cl_2p} is hyperbolic. 
Yet for mixed states \(\UU \in \cU \backslash (\PP_\phasel \cup \PP_\phaser)\) the system can become elliptic. 
The lack of global hyperbolicity is a major challenge, that prohibits the usage of standard  solvers. 
The domains \(\PP_\phasel\), \(\PP_\phaser\) form the state space for the two phases, which are marked by the subscripts \(\phasel\), \(\phaser\). 
Accordingly, we will speak of the \(\phasel\)-phase and the \(\phaser\)-phase.
We will use the abbreviated notation \(\phaserl\) to indicate both phases (separately).
 \subsection{The Liquid--Vapor Flow Model}
\label{sec:euler}

In this section we consider liquid--vapor flow as an explicit instance of the general two-phase flow model of Section~\ref{sec:generic_two_phase}.

\paragraph*{Bulk Equations}
We use the subscript \(\phasel\) to refer to the \emph{liquid phase}, and the subscript \(\phaser\) to refer to the \emph{vapor phase}. 
In order to model the flow in a domain \(\contdomain \subset \bR^d\) with a sharp interface \(\continterface(t)\), we follow the geometric setting described 
in Section~\ref{sec:generic_two_phase}.
We assume that the Euler equations hold in both bulk phases, i.e.
\begin{align} 
  \label{eq:full_euler}
 \begin{aligned}
  \partial_t \rho + \nabla \cdot (\rho \vv) &= 0, \\
  \partial_t (\rho \vv) + \nabla \cdot (\rho \vv \otimes \vv + p \mat{I}) &= \vec{0}, \\
  \partial_t E + \nabla \cdot ((E + p) \vv) &= 0,
 \end{aligned}
 \quad\text{ in } \contdomain_\phaserl(t) \text{ for } t \in (0,t_{\mathrm{end}}).  
\end{align}
Again, we make use of the abbreviated notation \eqref{eq:compact_phase_notation}. 
The unknowns of the system \eqref{eq:full_euler} are the fluid density \(\rho = \rho(\vx,t)\), the fluid momentum density \(\rho \vv = \rho\vv(\vx,t)\), the total energy density \(E = E(\vx,t)\), and the interface \(\continterface(t)\).
The total energy density satisfies \(E = \rho \varepsilon + \tfrac{1}{2} \rho \norm{\vv}^2\), with \(\varepsilon\) denoting the specific internal energy.
The \(d\)-dimensional identity matrix is written as \(\mat{I}\). 
The flux function \(\vf\) of \eqref{eq:full_euler}
is given by 
\begin{align}
  \label{eq:euler_flux}
  \vf(\UU) \coloneqq 
  \begin{pmatrix}
    \vm^\top \\ 
    \tfrac{1}{\rho} \vm \otimes \vm + p \mat{I} \\    
    (E + p) \vv^\top 
  \end{pmatrix},
\end{align}
with the state vector \(\UU \coloneqq (\rho, \vm, E)^\top\), where \(\vm \coloneqq \rho \vv\) denotes the momentum density.
\newline
To further close the system, the pressure \(p\), the specific internal energy \(\varepsilon\), and the temperature \(T\) are expressed in terms of thermodynamically independent quantities by an equation of state (\eosshort{}).
Many \eosshort{}  do not prescribe e.g. the pressure directly.  
Instead, the pressure can be expressed via  
\(p(\rho) = \rho^2 \tfrac{\dd}{\dd\!\rho}\psi(\rho)\), 
with \(\psi\) being  the specific Helmholtz free energy.
Likewise, the specific internal energy can be computed via \(\varepsilon(\rho, T) = \psi(\rho,T) - T \tfrac{\dd}{\dd\!T}\psi(\rho,T)\), as well as  the internal entropy \(S\) via \(S = S(\rho,T) = -\partial_T\psi(\rho,T)\).

As discussed in \cite{dafermos:hyperbolic:2016}, the system \eqref{eq:full_euler} is hyperbolic if the \eosshort{}  fulfills 
\begin{align}
  \label{eq:hyperbolicity_condition_euler}
  \begin{aligned}
    \partial_\rho p(\rho,S) &> 0, 
    &
    T(\rho,S) &> 0.   
  \end{aligned}
\end{align}
The admissibility criterion \eqref{eq:hyperbolicity_condition_euler} divides the state space \(\cU\), depending on the \eosshort{}, into admissible regions that are identified with different fluid phases. 
We denote the domain in the state space \(\cU\) that corresponds to the liquid phase by \(\PP_{\phasel}\), and for the vapor phase by \(\PP_{\phaser}\).

\begin{figure}[t]
\centering
\subcaptionbox{Pressure \(p = p(\rho, T)\)
\label{subfig:lj_eos:pressure}}[0.49\columnwidth]
{\includegraphics[width=0.45\columnwidth]{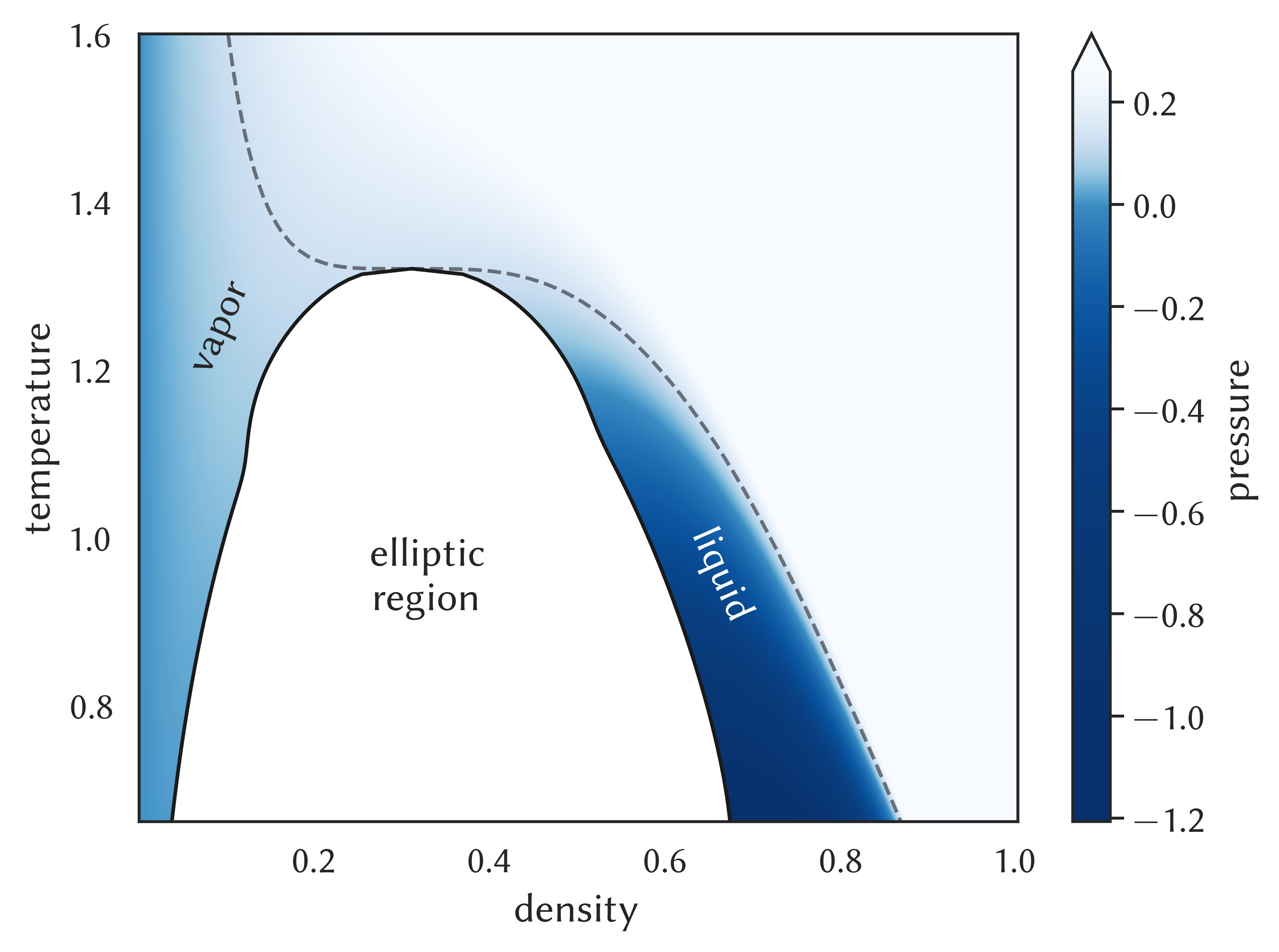}
}
\subcaptionbox{Temperature \(T = T(\rho, \varepsilon)\)
\label{subfig:lj_eos:temperature}}[0.49\columnwidth]
{\includegraphics[width=0.45\columnwidth]{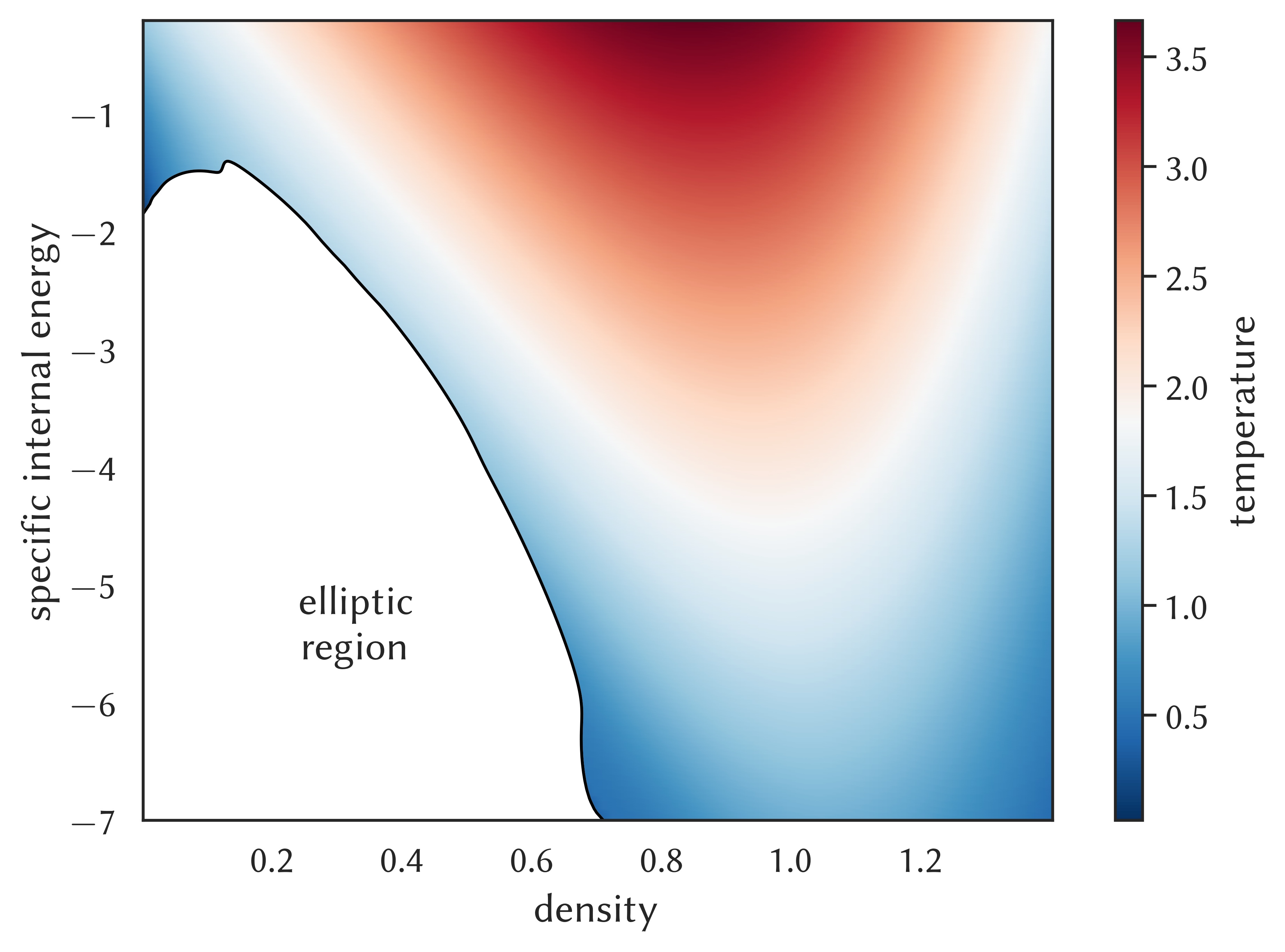}
}
\caption{Pressure and temperature of the Lennard--Jones \eosshort{}  \cite{thol.rutkai.ea:equation:2016}.
The dashed line represents the isobar of the critical pressure \(p_{\mathrm{crit}} = p(\rho_{\mathrm{crit}}, T_{\mathrm{crit}})\).
}
\label{fig:noniso_lennard_jones_eos}
\end{figure}

\paragraph*{Equation of State for a Lennard--Jones Fluid}
To describe the continuum-scale properties of the fluid in \eqref{eq:full_euler}, we have to specify an explicit \eosshort{}.
We choose the \eosshort{}  presented in \cite{thol.rutkai.ea:equation:2016} that is consistent with a Lennard--Jones fluid 
on the molecular scale, and can be used to describe e.g. noble gases.
A graphical representation of the pressure \(p = p(\rho,T)\) in terms of density \(\rho\) and temperature \(T\) is shown in the phase diagram in Figure~\ref{fig:noniso_lennard_jones_eos}.

\paragraph*{Interface Conditions}
The Rankine--Hugoniot conditions \eqref{eq:generic_rh} for the system \eqref{eq:full_euler} write as
\begin{align}
  \label{eq:euler_rankine_hugoniot}
\begin{aligned}
\jump[\big]{ \rho ( \vv \cdot \vn - s ) } &= 0, \\ 
\jump[\big]{ \rho ( \vv \cdot \vn - s ) \vv \cdot \vn + p } &= 0, \\
\jump[\big]{ E ( \vv \cdot \vn - s ) + p  \vv \cdot \vn } &= 0.
\end{aligned}
\end{align}
Note that we neglect interface curvature and surface tension effects. 
If these effects are taken into account, the right-hand sides of \eqref{eq:euler_rankine_hugoniot}\textsubscript{2,3} are non-zero and depend on the local curvature.

In contrast to hydrodynamical shock waves, phase boundaries usually occur as subsonic waves, see Figure~\ref{fig:generic_wave_pattern} for a typical Riemann wave pattern.
Hence, additional criteria have to be imposed, to ensure a unique solution \cite{abeyaratne.knowles:evolution:2006,truskinovsky:kinks:1993}.
In classical theories 
\cite{abeyaratne.knowles:kinetic:1991,lefloch.thanh:nonclassical:2000,lefloch.thanh:nonclassical:2001,lefloch.thanh:non:2002} this is done by prescribing kinetic relations like in \eqref{eq:generic_kinetic_relation}. 
For example, one can specify a driving force \(\cK \colon \bR \to \bR\) that describes the entropy production across the interface.
In this instance, the kinetic relation is given by
  \begin{align}
    \label{eq:euler_kinrel}
  \jump[\big]{ S } - \cK(j) = 0,   
 \end{align}  
with \(S = S(\rho,T)\) denoting the internal entropy of the fluid  and  \(j \coloneqq \rho_{\phaserl} (\vv_{\phaserl} \cdot \vn - s)\) the relative mass flux.

\section{The Microscale Interface Solver}
\label{sec:interface_solver}

The core problem of describing two-phase flow lies in the accurate description of the phase boundary dynamics. 
Within our framework, this can be realized by solving a planar two-phase Riemann problem, i.e. by solving \eqref{eq:generic_cl_2p} for the Riemann initial data 
\begin{align} \label{eq:riemanndata}
  \UU_0(x) = \begin{cases}
  \UU_{\phasel} & \text{ for } x < 0, \\
  \UU_{\phaser} & \text{ for } x > 0,
\end{cases}
\end{align}
with \(\UU_{\phaserl} \in \PP_\phaserl\) constant, and \(x = (\vx - \vec{\xi}) \cdot \vn\), for \(\vx \in \bR^d\), \(\vec{\xi} \in \continterface\).
Solving the two-phase Riemann problem analytically in the context of temperature-dependent liquid--vapor flow with sharp interfaces remains a largely unresolved task on the continuum-mechanical level. 
Despite the fact, that the two-phase Riemann problem has been widely analyzed in many contributions \cite{abeyaratne.knowles:kinetic:1991,yanagi:riemann:1992,colombo.priuli:characterization:2003,merkle.rohde:sharp:2007,hantke.dreyer.ea:exact:2013,thanh.vinh:riemann:2021}, 
there remain major open issues for realistic \eosplshort{} and temperature-dependent flow. 
One problem is that the physically relevant choice for the driving force \(\cK(j)\) in the kinetic relation \eqref{eq:euler_kinrel} is often unclear and depends on the type of fluid that is considered. 
Even more so, it is shown in \cite{thein:results:2018,hantke.thein:impossibility:2019} 
that there are phase boundaries that fulfill the basic Rankine--Hugoniot jump conditions but contradict entropy constraints in the first-order framework.
\newline 
These shortcomings of the continuum-scale model are overcome by our  multiscale approach.
In that we
utilize molecular dynamics simulations at the interface. 
Thus, the driving forces are inherently included on the particle scale and do not have to be modeled separately on the continuum-scale. 
Before discussing the specific molecular-scale interface solver in detail, we start by introducing the general concept of a microscale interface solver.

\subsection{The Microscale Interface Solver for the General Two-Phase Framework}

\paragraph*{Directional State Projections}
The core problem of the multiscale model-algorithm is to find solutions of the local Riemann problems along the interface \(\continterface(t)\). 
For this purpose, we regard the directional version of \eqref{eq:generic_cl_2p}. 
Here, we have to distinguish between conserved quantities in \(\UU\) that are scalar (e.g. mass), and others that are vector-valued (e.g. momentum). 
Assume that \(\UU = (u_1,\ldots, \allowbreak u_k,\allowbreak \vv_1, \ldots, \allowbreak \vv_l) \allowbreak \in \bR^m\) (after reordering), with \(u_i \in \bR\), \(i = 1, \ldots,k\), and \(\vv_i \in \bR^d\), \(i = 1, \ldots, l\), being primary unknowns of the system \eqref{eq:generic_cl_2p}. 
For \(\vn \in \bS^{d-1}\) the \emph{projected state in direction of \(\vn\)} and its \emph{perpendicular remainder} are defined by
\begin{align}
  \label{eq:directional_states}
  \begin{aligned}
  \UU_{\parallel\vn} &\coloneqq (u_1,\ldots, u_k, \vv_1 \cdot \vn, \ldots, \vv_l  \cdot \vn) \in \bR^{k+l}, \\
  \UU_{\perp\vn} &\coloneqq (0,\ldots, 0, \vv_1 - (\vv_1 \cdot \vn) \vn, \ldots, \vv_l  - (\vv_l \cdot \vn) \vn) \in \bR^{m}.
\end{aligned}
\end{align}
The \emph{directional state} \(\UU_{\parallel\vn}\) can be transferred back to its original form via the operator 
\begin{align}
  \label{eq:reverse_projection}  
    P_{\vn} \colon \bR^{k+l} \to \bR^m 
    : 
    (u_1,\ldots, u_k, v_1, \ldots, v_l)
    \mapsto     
    (u_1,\ldots, u_k, v_{1} \vn, \ldots, v_{l} \vn).  
\end{align}
Indeed, it holds \(\UU = P_{\vn}(\UU_{\parallel\vn}) + \UU_{\perp\vn}\).
For the projection in case of the  liquid--vapor flow model we refer to  \eqref{eq:noniso:multiscale:projected_state_variables}.

\paragraph*{Rotated Riemann Problem in Normal Direction \(\vn \in \bS^{d-1}\)}
Let \(x = (\vx - \vec{\xi}) \cdot \vn\), with \(\vx \in \bR^d\), \(\vec{\xi} \in \continterface\). 
The rotated initial value problem corresponding to \eqref{eq:generic_cl_2p} is given by 
\begin{align} \label{eq:rotated_cl}
  \begin{aligned}
    \partial_t \uu + \partial_x \vF(\uu; \vn) &= \vec{0}, 
    && \text{ in } \bR \times (0,t_{\mathrm{end}}), \\
    \uu(x,0) &= \uu_0(x,0),  
    && \text{ in } \bR,
  \end{aligned}
\end{align}
where \(\uu \coloneqq \UU_{\parallel\vn}\) is the  projected state in direction of \(\vn\).
For \(\UU_{\phasel}\), \(\UU_{\phaser} \in \cU\), we define the rotated Riemann initial data 
\begin{align} \label{eq:generic_rotated_riemanndata}
  \uu_0(x) = \begin{cases}
  \uu_{\phasel} & \text{ for } x < 0, \\
  \uu_{\phaser} & \text{ for } x > 0,
\end{cases}
\end{align}
where \(\uu_{\phasel} \coloneqq (\UU_{\phasel})_{\parallel\vn}\), \(\uu_{\phaser} \coloneqq (\UU_{\phaser})_{\parallel\vn}\)
are the projected states in direction of \(\vn\).

\begin{figure}[tb]
\centering
\includegraphics[width=0.33\columnwidth]{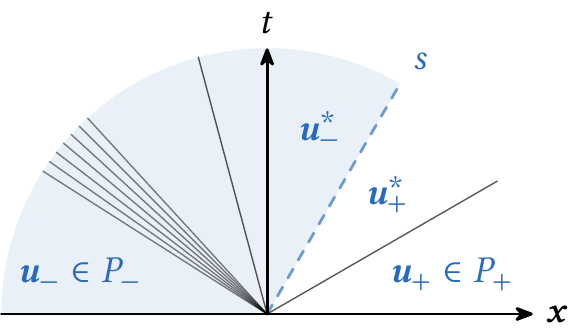}
\caption{
Sketch of a generic, one-dimensional Riemann wave pattern for liquid--vapor flow, resulting from the Riemann problem with initial data \(\uu_{\phasel} \in \pp_{\phasel}\) and \(\uu_{\phaser} \in \pp_{\phaser}\).
The dashed line marks the phase boundary, which is sharp as an additional discontinuous wave, and it forms the wave of interest.
}
\label{fig:generic_wave_pattern}
\end{figure}

\paragraph*{Wave of Interest}
On the continuum scale, the solution of the Riemann problem \eqref{eq:rotated_cl}, \eqref{eq:generic_rotated_riemanndata} consists of a composition of elementary waves, such as shocks, contacts, and rarefactions.
Out of them, we are interested in the dynamics of one specific, discontinuous wave of interest, which represents  the sharp interface.
We refer to Figure~\ref{fig:generic_wave_pattern} for an illustration of a typical wave pattern.
For the wave of interest we denote its two adjacent states by \(\uu^*_\phaser\), \(\uu^*_\phasel \in \bR^{k+l}\), and the wave speed by \(s \in \bR\).
This corresponds to a traveling wave
\begin{align} 
\label{eq:rotated_wave_of_interest_in_normal_direction}
\uu^*(x,t) =
\begin{cases}
   \uu^{*}_{\phasel} & \text{ for } x < st, \\
   \uu^{*}_{\phaser} & \text{ for } x > st.
\end{cases}
\end{align}
The traveling wave \eqref{eq:rotated_wave_of_interest_in_normal_direction} is a weak solution of the rotated Riemann problem with the initial states \(\uu^*_\phaser\), \(\uu^*_\phasel\).
As a prerequisite for weak solutions, the Rankine--Hugoniot conditions 
\begin{align}
  \label{eq:rankine_hugoniot_condition}
 \vF(\uu^{*}_{\phasel}; \vn(\vec{\xi})) - \vF(\uu^{*}_{\phaser}; \vn(\vec{\xi}))
 = s\,(\uu^{*}_{\phasel} - \uu^{*}_{\phaser}),
\end{align}
have to be fulfilled at the interface \(\continterface(t)\).

Projected back to the full-dimensional system, the wave of interest at the interface \(\continterface(t)\) consists of the states 
\(\UU^{*}_{\phaserl} \coloneqq P_{\vn}(\uu^*_\phaserl) + (\UU_\phaserl)_{\perp\vn} \in \PP_\phaserl\), 
and the wave speed \(s \in \bR\).
The tuple \((\UU^{*}_{\phasel}, \UU^{*}_{\phaser}, s)\) forms a traveling wave
\begin{align} 
\label{eq:wave_of_interest_in_normal_direction}
\UU^*(\vx,t; \vec{\xi}) =
\begin{cases}
           \UU^{*}_{\phasel} & \text{ for } \vx \cdot \vn(\vec{\xi}) < st, \\
           \UU^{*}_{\phaser} & \text{ for } \vx \cdot \vn(\vec{\xi}) > st,
          \end{cases}
\end{align}
which is itself a weak solution of \eqref{eq:generic_cl_2p} with discontinuous initial data consisting of the two adjacent states.

\paragraph*{General Interface Solver}
Our main interest lies in the construction of an interface solver, that computes the adjacent states at the interface, and the interface speed, for given initial data \eqref{eq:generic_rotated_riemanndata}. 
We define such a solver in the following way.

\begin{definition}[General Interface Solver]
\label{def:generic_interface_solver}
Let \(\UU_\phasel \in \PP_\phasel\), \(\UU_\phaser \in \PP_\phaser\) and \(\vn \in \bS^{d-1}\) be given. 
We assume that for the corresponding Riemann initial data \eqref{eq:generic_rotated_riemanndata} a unique weak solution of the rotated system \eqref{eq:rotated_cl} exists, and is composed of elementary waves that include a (discontinuous) wave of interest, defined by the adjacent states \(\uu^*_\phaser, \uu^*_\phasel \in \bR^{k+l}\) and wave speed \(s \in \bR\). 
Then, if the states \(\UU^{*}_{\phaserl} \coloneqq P_{\vn}(\uu^*_\phaserl) + (\UU_\phaserl)_{\perp\vn}\) 
are elements of \(\PP_\phaserl\), 
we call the mapping 
\begin{align}
  \genericsolversym \colon \PP_\phasel \times \PP_\phaser \times \bS^{d-1} \to \PP_\phasel \times \PP_\phaser \times \bR 
  :
 (\UU_\phasel, \UU_\phaser, \vn) \mapsto (\UU^{*}_{\phasel}, \UU^{*}_{\phaser}, s)
\end{align}
a \emph{general interface solver} of the system \eqref{eq:generic_cl_2p}.
\end{definition}

As discussed in Section~\ref{sec:euler}, such an interface solver is not easily available in a closed formulation as an exact solution of the Riemann problem.

\begin{remark}
\label{remark:interface_solver}
For the construction of the interface solver, we are mostly interested in the dynamics of the Riemann problem \eqref{eq:generic_rotated_riemanndata} of the rotated system \eqref{eq:rotated_cl}. 
If the flux function \(\vf\) is rotationally invariant, or \(d=1\), it suffices to construct an interface solver of the form 
\begin{align}
  \label{eq:riemannsolver:rotated}
  \rmicro \colon
  \pp_\phasel \times \pp_\phaser
  \to \pp_\phasel \times \pp_\phaser \times \bR
  :
  (\uu_{\phasel}, \uu_{\phaser}) 
  = (\uu^{*}_{\phasel}  \uu^{*}_{\phaser}, s), 
\end{align}
with \(\pp_{\phaserl}\) being the corresponding projections of the phase state spaces \(\PP_{\phaserl}\).  
\end{remark}
 \subsection{The Molecular-Scale Interface Solver}
\label{sec:md}

One of the most detailed ways to describe fluid flow on a microscopic scale is to perform 
\mdshort{}
simulations, where each molecule is represented as a single particle. 
This has already been successfully applied to investigate liquid--vapor flow \cite{matsumoto:molecular:1998,homes.heinen.ea:evaporation:2020}, albeit, mainly focused on equilibrium computations.
\newline
In this section, first, we summarize some basics for \mdshort{} simulations (see also \cite{allen.tildesley:computer:2017} for a general introduction).
For details on our specific \mdshort{} simulations we refer to Appendix~\ref{appendix:md_elements}.
Then, we apply \mdshort{} simulations to construct a novel microscale interface solver that describes the non-equilibrium dynamics of the phase boundary. 

\paragraph*{\mdshort{} System}
For an \mdshort{} particle system we consider \(\nparticles \in \bN\) particles with masses \(m_i > 0\), 
positions \(\vx_i = (x_i, y_i, z_i) \in \bR^3\), 
velocities \(\vv_i \in \bR^3\), 
and accelerations \(\va_i \in \bR^3\), with \(i = 1, \ldots, \nparticles\), 
inside a molecular-scale domain \(\mddomain \subset \bR^3\). 
The particles interact via an inter-particle, two-body potential \(\phi \colon \allowbreak (0,\infty) \to \bR\). 
We consider the Lennard--Jones potential 
\begin{align} 
  \label{eq:md:lennard_jones_potential}
  \phi(r) = 4 \varepsilon \left( \left( \frac{\sigma}{r} \right)^{12} - \left( \frac{\sigma}{r} \right)^6 \right),
\end{align}
where \(r > 0\) denotes the distance between two particles. 
The parameters  \(\sigma, \varepsilon > 0\) specify the type of particles and therefore the type of fluid. 
For \(\sigma = 1\) and \(\varepsilon = 1\) we get a Lennard--Jones fluid whose properties are  described by the \eosshort{} \cite{thol.rutkai.ea:equation:2016} of Section \ref{sec:euler}.
The force between two arbitrary particles \(i \neq j\) computes as
\begin{align*}
 \vf_{ij} = - \nabla_{\vx_i} \phi\bigl(\norm{\vx_i - \vx_j}_2\bigr).
\end{align*}
With this, we can formulate the equations of motion for each particle \(i = 1,\ldots,\nparticles\) and time \(\mdtimevar \geq 0\) by
\begin{align} 
  \label{eq:md:equation_of_motion}
\begin{aligned}
 \tfrac{\dd }{\dd \mdtimevar} \vx_i(\mdtimevar)  &= \vv_i(\mdtimevar), 
 & \tfrac{\dd }{\dd \mdtimevar} \vv_i(\mdtimevar) &= \va_i(\mdtimevar), 
  & \va_i(\mdtimevar) &= \sum_{j \neq i} \vf_{ij}(\mdtimevar) / m_i,
  \end{aligned}
\end{align}
subject to the initial conditions
\begin{align} 
  \label{eq:md:equation_of_motion:initial}
\begin{aligned}
 \vx_i(0)  &= \vx^{0}_{i}, 
 & 
  \vv_i(0) &= \vv^{0}_{i}.
  \end{aligned}
\end{align}

Note that the \mdshort{} simulations are always three-dimensional, as only in this case the micro-scale fluid is consistent with the continuum-scale fluid and its \eosshort{}.
The discretization of the particle system \eqref{eq:md:equation_of_motion} is given in Algorithm~\ref{alg:md_step}. 
By using the CUDA framework \cite{nickolls.buck.ea:scalable:2008} we are able to utilize graphics processing units (GPUs) that achieve a high computational performance for such tasks.

\begin{figure}[tp]
\centering
\includegraphics[width=0.8\columnwidth]{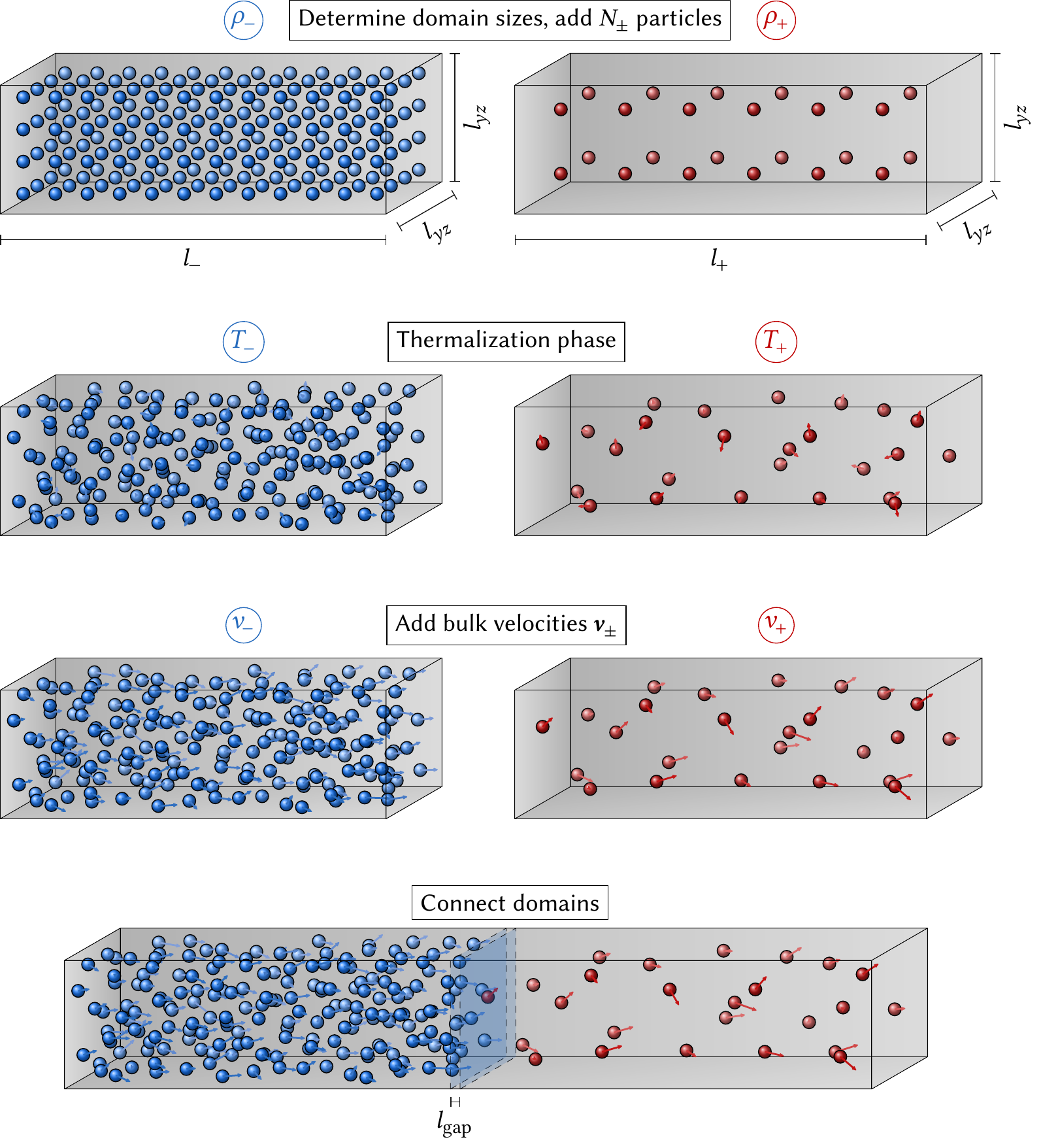}
\caption{Sketch of the initialization of a molecular-scale Riemann problem.}
\label{fig:md_riemann_initialization}
\end{figure}

\paragraph*{Molecular-Scale Riemann Problem}
Our main interest lies in computing Riemann solutions on the molecular scale for the continuum-scale initial conditions \(\uu_{\phasel} = (\rho_{\phasel}, v_{\phasel}, T_{\phasel})\) in the liquid phase and \(\uu_{\phaser} = (\rho_{\phaser}, v_{\phaser}, T_{\phaser})\) in the vapor phase, as in \eqref{eq:generic_rotated_riemanndata}.
Both phases are divided by a molecular-scale interface which we consider as a plane positioned at \(\Gamma(0) \in \bR\) in \(x\)-direction.
\newline
To set up such a molecular-scale Riemann problem, we choose a fixed number of particles \(\nparticles \in \bN\) and distribute them into two subdomains with fixed width in transversal directions, such that the desired density is achieved and the ratio of the subdomain lengths in \(x\)-direction is fixed.
\newline
In the next step, both particle systems have to undergo a thermalization phase, to attain an appropriate thermodynamic equilibrium state corresponding to \(\uu_{\phaserl}\). 
In this phase, we consider each subsystem separately with periodic boundary conditions. 
Then, we run the \mdshort{} simulation for a fixed number of time steps, while applying a thermostat in each subdomain, to attain the desired temperatures \(T_{\phaserl}\). 
Afterwards, the bulk velocities \(v_{\phaserl}\) are added to all particle velocities in each respective subdomain. 
This realizes the continuum-scale fluid velocities on the particle scale.
Finally,
both domains are combined by putting them next to another, while retaining a small gap \(l_{\mathrm{gap}} = 2^{\frac{1}{6}} \sigma \) in between, in order to avoid overlapping particles. 
In that way, we have set up the initial particle system. 
The whole initialization procedure is
described in Algorithm~\ref{alg:md_riemann_init} in Appendix~\ref{appendix:md_elements} --- 
for an illustration we refer to Figure~\ref{fig:md_riemann_initialization}.

\paragraph*{Interface Tracking}
\label{sec:md_tracking}
To get the position of the liquid--vapor interface, we exploit the fact that there is a large density jump from the liquid to the vapor phase. 
The first step is to transfer the particle distribution \((\vx_i(\mdtimevar))_{i=1,\ldots,\nparticles}\) to a continuous density distribution \(\widehat{\rho}(x, \mdtimevar)\) along the \(x\)-axis by employing 
radial basis functions (RBFs)
\begin{align} 
  \label{eq:md_density_rbf}
  \widehat{\rho}(x,\mdtimevar) 
  \coloneqq 
  \sqrt{ \frac{\gamma_{\mathrm{RBF}} }{\pi} }
  \sum_{i=0}^{\nparticles} m_i \exp\bigl(-\gamma_{\mathrm{RBF}} (x - x_i(\mdtimevar))^2  \bigr),
\end{align}
where \(\gamma_{\mathrm{RBF}} > 0\) is the RBF-parameter, controlling the smoothing.
The derivative with respect to \(x\) computes as
\begin{align} 
  \label{eq:md_density_rbf_dx}
  \partial_x \widehat{\rho}(x,\mdtimevar) 
  = 
  - 2 \sqrt{ \frac{\gamma_{\mathrm{RBF}}^3}{\pi} }
  \sum_{i=0}^{\nparticles} m_i (x - x_i(\mdtimevar)) \exp\bigl(-\gamma_{\mathrm{RBF}} (x - x_i(\mdtimevar))^2  \bigr).
\end{align}
During the \mdshort{} simulation (computed using the time-step Algorithm~\ref{alg:md_step}) 
we compute the interface position at every time step \(\mdtimevar_n\),
by finding the maximum of \(\abs{\partial_x \widehat{\rho}(x,\mdtimevar_n)}\), 
yielding the new interface position \(\mdinterface(\mdtimevar_n) \in \bR\). 
For the local optimization we use basic optimization routines implemented in e.g. \cite{king:dlib:2009}.

After obtaining the interface position \(\mdinterface(\mdtimevar_n)\), we infer the interface speed by performing a simple linear regression over the last \(n_{\mathrm{buffer}} \in \bN\) interface positions \((\mdtimevar_j, \mdinterface(\mdtimevar_j))\), \(j = \max(n - n_{\mathrm{buffer}}, 0), \ldots, n\). 
The slope of the fitted line provides the approximation for the interface speed \(s = s(\mdtimevar_n)\) at time \(\mdtimevar_n\).
The whole procedure is described in Algorithm~\ref{alg:md_tracking} in Appendix~\ref{appendix:md_elements}.

\paragraph*{Interface State Sampling}
For the multiscale model, it is necessary to compute the fluid states \(\uu^*_\phaserl = (\rho^*_\phaserl, \vm^*_\phaserl, T^*_\phaserl)\) adjacent to the interface, see Figure~\ref{fig:generic_wave_pattern}.
To compute these states during the \mdshort{} simulation, we define a sampling region \(\Sigma_{\phaserl}(\mdtimevar)\) near the interface for each phase:
\begin{align}
  \Sigma_\phasel(\mdtimevar) & \coloneqq [\mdinterface(\mdtimevar) - w_{\mathrm{sr}} - o_{\mathrm{sr}}, \mdinterface(\mdtimevar) - o_{\mathrm{sr}}], 
  &
  \Sigma_\phaser(\mdtimevar) & \coloneqq [\mdinterface(\mdtimevar) + o_{\mathrm{sr}}, \mdinterface(\mdtimevar) + w_{\mathrm{sr}} + o_{\mathrm{sr}}],   
\end{align}
where \(w_{\mathrm{sr}} > 0\) denotes the width of the sampling region, and \(o_{\mathrm{sr}} \geq 0\) the interface offset. 
The width of the sampling region is determined heuristically, but it has a large impact on the simulation results. 
Making it too small results in high statistical errors as the number of particles in the sampling domain decreases.
In contrast, making it too large may give incorrect wave states if other waves are included. 
This means the sampling region should not be larger than the wave state plateaus at the interface.
\newline
The interface offset accounts for the small transition region near the interface, which should be avoided in the averaging procedure.

\begin{figure}[t]
\centering
\includegraphics[width=0.75\columnwidth]{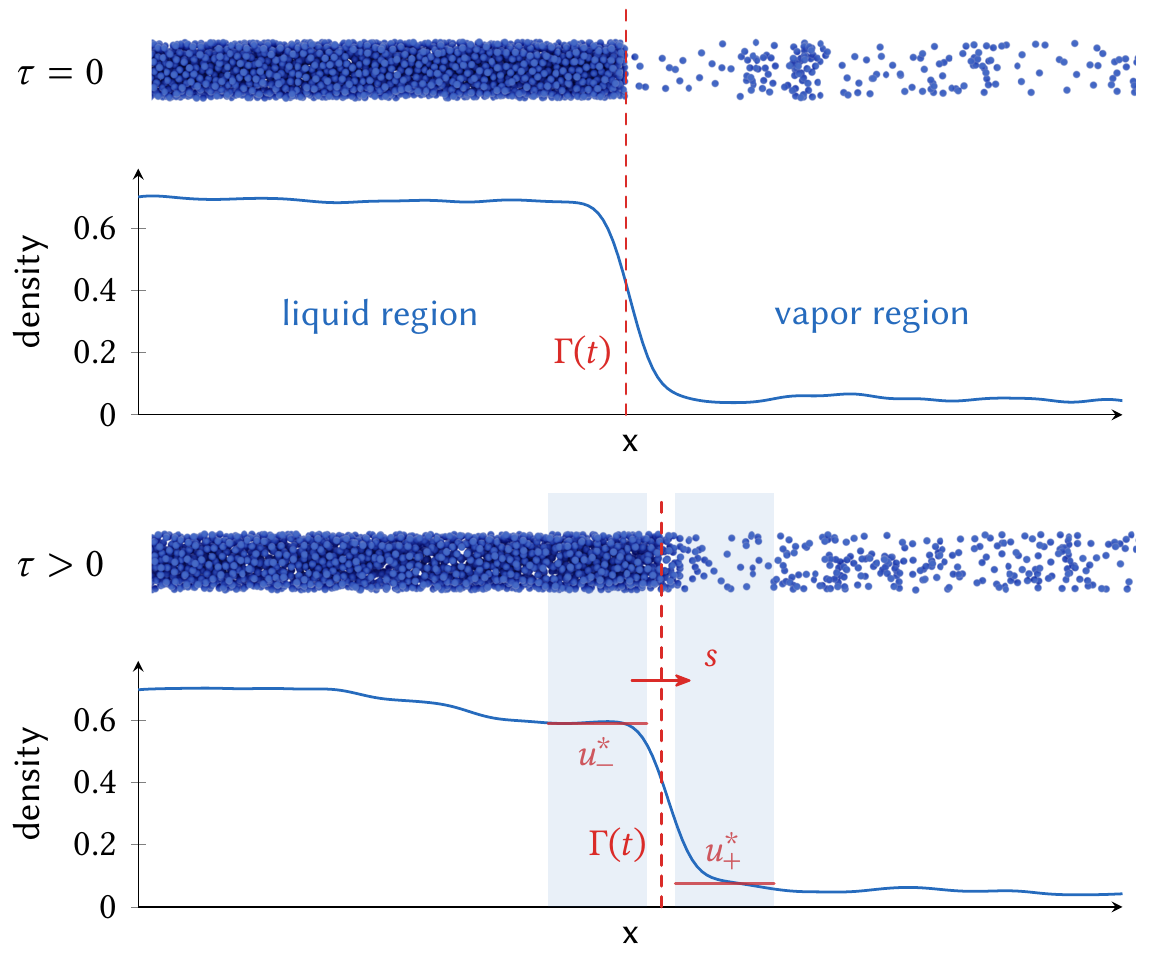}
\caption{Sampling of the density interface wave states during an \mdshort{} simulation.
The shaded regions mark the sampling domains \(\Sigma_\phaserl(\mdtimevar)\).
}
\label{fig:md_density_profile}
\end{figure}

Having defined the sampling regions, we can compute the interface wave states by employing the Irving--Kirkwood formulas \cite{irving.kirkwood:statistical:1950,e:principles:2011} 
\begin{align} 
  \label{eq:md_interface_states}
\begin{aligned}  
  \rho^*_\phaserl(\mdtimevar)
  & = \operatorname{vol}(\Sigma_{\phaserl})^{-1} \sum_{x_i(\mdtimevar) \in \Sigma_{\phaserl}} m_i ,
  \\
  \vm^*_\phaserl(\mdtimevar)
  & = \operatorname{vol}(\Sigma_{\phaserl})^{-1} \sum_{x_i(\mdtimevar) \in \Sigma_{\phaserl}} m_i \vv_i(\mdtimevar), 
  \\
  \vv^*_\phaserl(\mdtimevar)
& = \Bigl( \sum_{x_i(\mdtimevar) \in \Sigma_{\phaserl}} m_i \Bigr)^{-1} \sum_{x_i(\mdtimevar) \in \Sigma_{\phaserl}} m_i \vv_i(\mdtimevar), 
  \\
  T^*_\phaserl(\mdtimevar)
  & = \abs{x_i(\mdtimevar) \in \Sigma_{\phaserl}}^{-1} \frac{1}{d} \sum_{x_i(\mdtimevar) \in \Sigma_{\phaserl}} m_i (\vv_i(\mdtimevar) - \vv^*_\phaserl(\mdtimevar)) \cdot (\vv_i(\mdtimevar) - \vv^*_\phaserl(\mdtimevar)),  
\end{aligned}
\end{align}
where \(\operatorname{vol}(\Sigma_{\phaserl}) = w_{\mathrm{sr}} h_{\mathrm{dom}}^2\) denotes the volume of the sampling region, and \(\abs{x_i(\mdtimevar) \in \Sigma_{\phaserl}}\) the number of particles inside it. 
In that way, we obtain the local interface states at every point in time during a \mdshort{} simulation. 
To get even more reliable results, we perform time averaging over a fraction \(\alpha_{\mdtimevar\mathrm{-smpl}} \in [0,1]\) of the total simulation time \(\mdtimeend > 0\).
\newline 
A graphical illustration of the interface state sampling and the interface tracking procedure is shown in Figure~\ref{fig:md_density_profile}.

\paragraph*{Molecular-Scale Interface Solver}
By applying the methods described in the previous sections, we are able to initialize a molecular-scale Riemann problem, run \mdshort{} simulations, and compute the interface dynamics. 
The whole procedure results in an interface solver (of the same form as  \eqref{eq:riemannsolver:rotated})
\begin{align} \label{eq:md_noniso_riemann_solver}
  \rmd(\rho_\phasel, v_\phasel, T_\phasel, \rho_\phaser, v_\phaser, T_\phaser) 
  &= 
  (\rho^*_\phasel, v^*_\phasel, T^*_\phasel, \rho^*_\phaser, v^*_\phaser, T^*_\phaser, s).
\end{align}
It maps the Riemann initial data to the interface wave states and the interface speed. 
The algorithmic description of \eqref{eq:md_noniso_riemann_solver} follows.

\begin{boxalgorithm}[label=alg:md_riemann]{Molecular-Scale Interface Solver}
\textbf{Input:}
  continuum fluid states \((\rho_\phasel, v_\phasel, T_\phasel)\), \((\rho_\phaser, v_\phaser, T_\phaser)\).

\textbf{Parameters:}
  total number of particles \(\nparticles\),
  processing frequency \(\frequencyvar_{\mathrm{pr}} \in \bN\),
number of time steps \(\mdstepend\),
  time step \(\Delta \mdtimevar > 0\),
  time sampling ratio \(\alpha_{\mdtimevar\mathrm{-smpl}} \in [0,1]\),
  sampling region width \(w_{\mathrm{sr}} > 0\),
  sampling region interface offset \(o_{\mathrm{sr}} \geq 0\),
  thermostat region offset ratio \(\alpha_{\mathrm{therm}} \in [0,1]\).
\tcbsubtitle{Algorithm}
\begin{itemize}  
 \item Set up initial particle configuration of the molecular-scale Riemann problem according to the continuum fluid states \((\rho_\phasel, v_\phasel, T_\phasel)\), \((\rho_\phaser, v_\phaser, T_\phaser)\) (see Algorithm~\ref{alg:md_riemann_init}). 
\item For \(n = 1,\ldots, \mdstepend\):
\begin{itemize}
  \item Compute \(\frequencyvar_{\mathrm{pr}}\) new time steps according to Algorithm~\ref{alg:md_step}.   
  \item Compute new interface position \(\mdinterface(\mdtimevar_n)\) and interface speed \(s(\mdtimevar_n)\) according to Algorithm~\ref{alg:md_tracking}.
  \item Compute interface states \((\rho^*_\phaserl, \vv^*_\phaserl, T^*_\phaserl)(\mdtimevar_n)\) as described in \eqref{eq:md_interface_states}.
\end{itemize}
 \item Compute the time averaged values \((\rho^*_\phaserl, \vv^*_\phaserl, T^*_\phaserl, s)\) from \((\rho^*_\phaserl, \vv^*_\phaserl, T^*_\phaserl)(\mdtimevar_n)\), \(s(\mdtimevar_n)\) for all \(n \in \{\widetilde{\alpha}_{\mdtimevar\mathrm{-smpl}} \mdstepend, \ldots, \mdstepend\}\). 
\end{itemize}
\textbf{Result:} interface states \((\rho^*_\phasel, \vv^*_\phasel, T^*_\phasel)\), \((\rho^*_\phaser, \vv^*_\phaser, T^*_\phaser)\), and interface speed \(s\).
\end{boxalgorithm}

\section{Surrogate Microscale Interface Solver}
\label{sec:surrogate}
In most cases, microscale simulations or the computation of the interface dynamics in general are computationally expensive. 
To address this problem, we employ \emph{surrogate solvers} that rely on machine learning algorithms, such as neural networks. 
The surrogate solver 
is denoted by 
\begin{align} \label{eq:generic_surrogate}
\networkF_{\nnparams} \colon \pp_{\phasel} \times \pp_{\phaser} \to \pp_{\phasel} \times \pp_{\phaser} \times \bR: (\uu_{\phasel}, \uu_{\phaser}) \mapsto (\uu^*_{\phasel}, \uu^*_{\phaser}, s),
\end{align}
where \(\nnparams\) constitutes the parameters of the surrogate model, which are determined by the training\slash{}optimization procedure.
For neural networks these parameters consist of the weights and biases of the network.   
\newline 
By using such a surrogate solver we can reduce the computational time from a couple of minutes for a full microscale simulation to a fraction of milliseconds.
This comes at the costs of generating a data set and training the surrogate solver. 
Nonetheless, applying surrogate solvers usually pays off after several evaluations, compared to more expensive microscale interface solvers \(\rmicro\).
For a more detailed discussion of the computational gains we refer to 
Section~\ref{sec:results:noniso:3d}.

\paragraph*{Constraint-Aware Neural Networks}
Most traditional methods for generating a surrogate solver \(\networkF_{\nnparams}\) do not incorporate physical constraints, such as the Rankine--Hugoniot conditions \eqref{eq:rankine_hugoniot_condition}
at the interface. 
This is primarily due to the fact that most methods merely focus on fitting the data, without trying to incorporate underlying physical constraints. 
Additionally, the microscale solver itself might not fulfill such constraints. 
This may be due to dissipative effects that appear on the microscale (but not in the continuum-scale model). 
Furthermore, statistical noise on the output values can invalidate constraints, which is the case for sampling from \mdshort{} simulations.
\newline
We mitigate this problem by employing \emph{constraint-aware neural networks}, see \cite{magiera.ray.ea:constraint:2020},
as our surrogate solver.
They are able to resolve physical constraints exactly by adding a constraint-resolving layer to a neural network, which implements a constraint-resolving function \(\Psi(\datz) = \daty\). 
To this end, we define the output variables \(\datz\)  of the previous part of the neural network, i.e. \(\networkG_{\nnparams}(\datx) = \datz\), and the mapping \(\Psi(\datz) = \daty\), such that  \(\daty\) fulfills the physical constraints.
The explicit form of the constraint-resolving function \(\Psi\) depends on the specific model. 
A particular choice for \(\Psi\) is given in \eqref{eq:constraint_resolving_function:euler} for liquid--vapor flow.

\paragraph*{Surrogate Solver Generation}
In this section, the generation of the surrogate solver is outlined. 
For an in-depth description of the employed machine learning algorithm we refer to \cite{magiera.ray.ea:constraint:2020}.

\paragraph*{Data Set Generation}
To generate a surrogate solver, we need a data set 
\(\datasetvar \coloneqq \{(\datx_i, \daty_i) \setsep  i = \allowbreak 1, \allowbreak \ldots, \allowbreak N_{\mathrm{data}}\}\) 
describing the input--response-relation of the microscale interface solver we intend to substitute. 
Here, the input 
\(\datx_i \coloneqq (\uu_{\phasel}, \uu_{\phaser})_i\) 
is a pair of Riemann initial data, and 
\(\daty_i \coloneqq (\uu^*_{\phasel}, \uu^*_{\phaser}, s)_i\) 
constitutes the corresponding output of the microscale interface solver \(\rmicro\), i.e. \(\daty_i = \rmicro(\datx_i)\).
Therefore, the first step is to generate an input data set 
\(\datasetvar_{\mathrm{in}} = \{ (\uu_{\phasel}, \uu_{\phaser})_i  \setsep   i = \allowbreak 1, \allowbreak \ldots, \allowbreak N_{\mathrm{data}} \}\).
Without active\slash{}online sampling methods (see Remark~\ref{remark:active_sampling}), this step involves choosing a sampling domain 
\(B_{\mathrm{in}} \coloneqq B_{\phasel} \times B_{\phaser} \subset \pp_{\phasel} \times \pp_{\phaser}\) 
that serves as the bounding domain for the input data, which means that all input points \(\datx_i\) will lie in \(B_{\mathrm{in}}\).
This is a crucial step, as the choice dictates for what kind of flow regimes the surrogate solver is valid. 
After determining such a sampling domain, input points for the microscale simulations are sampled in \(B_{\mathrm{in}}\) 
using the sampling Algorithm~\ref{alg:poisson_disc}.
This results in the input data set \(\datasetvar_{\mathrm{in}}\).
\newline
Afterwards, each input point \(\datx_i = (\uu_{\phasel}, \uu_{\phaser})_i \in \datasetvar_{\mathrm{in}}\) is evaluated by the microscale interface solver, i.e.
\begin{align*}
 \rmicro(\datx_i) = \daty_i = (\uu^*_{\phasel}, \uu^*_{\phaser}, s)_i,
\end{align*}
resulting in the output states \(\uu^*_{\phasel}\), \(\uu^*_{\phaser}\),  and the interface speed \(s \in \bR\).
This yields the corresponding output data set \(\datasetvar_{\mathrm{out}} = \{ (\uu^*_{\phasel}, \uu^*_{\phaser}, s)_i \setsep  i = \allowbreak 1, \ldots, \allowbreak N_{\mathrm{data}} \}\). 
Combining both data sets yields the complete data set \(\datasetvar \coloneqq \{(\datx_i, \daty_i) \setsep   i = \allowbreak 1, \ldots, \allowbreak N_{\mathrm{data}}\}\) that describes the input--response-relation of the microscale interface solver.

A large enough number of samples \(N_{\mathrm{data}}\) for the data set \(\datasetvar\) is necessary in order to apply the surrogate solver successfully. 
Note, that in this work we choose the data set size \(N_{\mathrm{data}}\) purely from an empirical standpoint.
Of course, the size is not the only quality measure of a data set, the distribution of the input variables \(\datx_i\) and inherent noise of the output variables \(\daty_i\) play a role just as important. 
In order to achieve a well-spaced input distribution, 
we apply the sampling Algorithm~\ref{alg:poisson_disc}.
The output noise can be reduced by performing several experiments for each input sample point \(\datx_i\) and averaging the outputs, resulting in an output label \(\daty_i\).

\subsection{Surrogate Solver for Liquid--Vapor Flow}
\label{sec:surrogate:liquid-vapor}
Instead of running the whole \mdshort{} simulation (Algorithm~\ref{alg:md_riemann}) for evaluating \(\rmd\) in \eqref{eq:riemannsolver:rotated}
we employ a surrogate solver. 
To decrease the dimensionality of the input space of the surrogate solver, we 
 exploit the fact that as a particle system  \mdshort{} are invariant with respect to a moving reference frame. 
 Thus, we choose the \(\phasel\)-phase velocity \(v_{\phasel} = \frac{m_\phasel}{\rho_\phasel}\) as the reference velocity \(v_{\mathrm{ref}} \coloneqq v_\phasel \) and consider the fluid velocities
\(\widetilde{v}_\phaserl = v_{\phaserl}  - v_{\mathrm{ref}}\).
Consequently, the surrogate solver takes the form 
\begin{align}  
  \networkF_{\nnparams}(\datx) 
  = \networkF_{\nnparams}(\rho_\phasel, T_{\phasel}, \rho_\phaser, v_\phaser, T_{\phaser})  
\approx 
  \rmd(\rho_\phasel, 0, T_{\phasel}, \rho_\phaser, v_\phaser, T_{\phaser})  
= 
  (\rho^*_\phasel, m^*_\phasel, T^*_\phasel, \rho^*_\phaser, m^*_\phaser, T^*_\phaser).
\end{align}
The underlying data set, as well as some remarks regarding the training of the network are found in Appendix~\ref{appendix:surrogate_details_euler}.

\begin{figure}[tb]
\centering
\includegraphics[width=0.5\columnwidth]{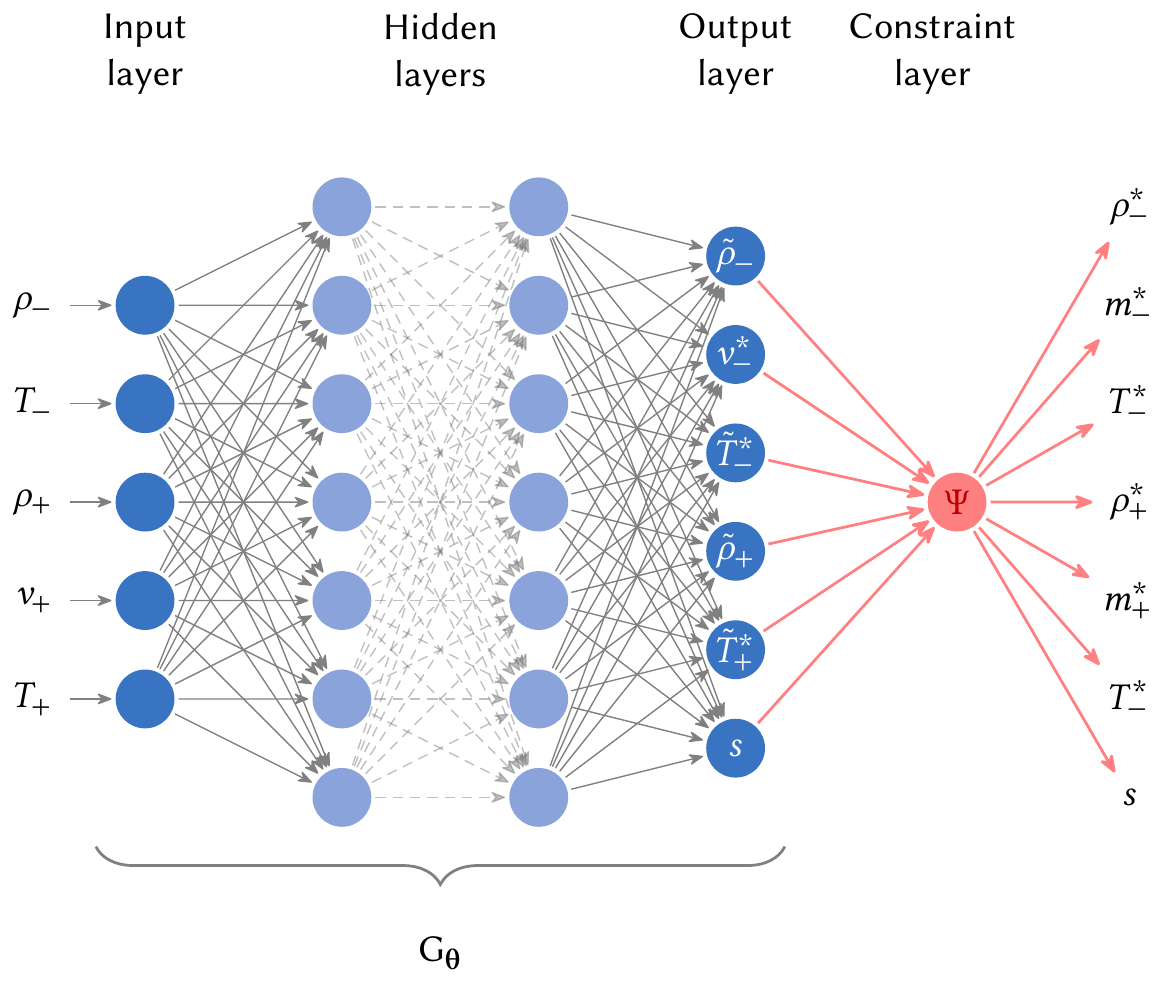}
\caption{
Graphical representation of the \cres-network for the liquid--vapor flow model. 
}
\label{fig:cres_noniso}
\end{figure}

As discussed in Section~\ref{sec:surrogate}, generic surrogate solvers \(\networkF_{\nnparams}\) are not tailored to uphold the Rankine--Hugoniot condition \eqref{eq:euler_rankine_hugoniot}, which include e.g. mass balance across the phase boundary. 
We aim to include mass conservation and positivity of the density in the surrogate solver \(\networkF_{\nnparams}\). 
For this purpose, we employ \cres-networks \cite{magiera.ray.ea:constraint:2020}
as our surrogate solver.
To build the necessary constraint-resolving function \(\Psi\), consider at first the mass conservation across the interface 
 \eqref{eq:euler_rankine_hugoniot}\textsubscript{1} rewritten as
\begin{align}
  \label{eq:mass_conservation:euler}
  m_\phaser =  m_\phasel - s (\rho_\phasel - \rho_\phaser).
\end{align}
This identity is used in the constraint-resolving function \(\Psi\) to account for the mass conservation across the interface. 
Furthermore, we ensure that the returned densities 
remain positive. 
\newline
The \cres-network consists of an intermediate network \(\networkG_{\nnparams}(\datx) = \datz\) with the output 
\begin{align}
  \datz \coloneqq ( \widetilde{\rho}_\phasel, v^{*}_{\phasel}, \widetilde{T}^{*}_\phasel, \widetilde{\rho}_\phaser, \widetilde{T}^{*}_\phaser, s ),
\end{align} 
where \(\widetilde{\rho}_\phaserl\) are the intermediate densities, 
\(\widetilde{T}^{*}_\phaserl\) the intermediate phase temperatures, 
\(v^{*}_{\phasel}\) the fluid velocity of the \(\phasel\)-phase, and \(s\) the interface speed.
Then, we define the function \(\Psi\) as
\begin{align}
  \label{eq:constraint_resolving_function:euler}
  \begin{aligned}
  \MoveEqLeft
  \Psi ( \widetilde{\rho}_\phasel, v^{*}_{\phasel}, \widetilde{T}^{*}_\phasel, \widetilde{\rho}_\phaser, \widetilde{T}^{*}_\phaser, s ) 
  \\
  &\coloneqq 
  \bigl( \abs{\widetilde{\rho}_\phasel}, \abs{\widetilde{\rho}_\phasel} v^{*}_{\phasel}, 
  \abs{\widetilde{T}^{*}_\phasel},
  \abs{\widetilde{\rho}_\phaser},
  \abs{\widetilde{\rho}_\phasel} v^{*}_{\phasel} - s(\abs{\widetilde{\rho}_\phasel} - \abs{\widetilde{\rho}_\phaser}),
  \abs{\widetilde{T}^{*}_\phaser}, 
  s
  \bigr) \\
  &= ( \rho^{*}_\phasel, m^{*}_\phasel,  T^{*}_\phasel, \rho^{*}_\phaser, m^{*}_\phaser, T^{*}_\phaser, s ).
\end{aligned}
\end{align}
This choice ensures that the output densities \(\rho^*_\phaserl\) and temperatures \(\widetilde{T}^{*}_\phaserl\) are positive, and by means of \eqref{eq:mass_conservation:euler} the mass across the interface is conserved. 
The general structure of the corresponding \cres-network is shown in Figure~\ref{fig:cres_noniso}.

In the end, a single evaluation of the complete surrogate solver takes less than a millisecond
whereas a single evaluation of \(\rmd\) takes between \num{9} and \SI{11}{\minute}. 
A more detailed discussion of the computational times is found in Section~\ref{sec:results:noniso:3d}.

\section{Finite Volume Method on Interface-Preserving Moving Meshes}
\label{sec:ipvm_algorithm}

\begin{figure}[tbp]
\centering
\includegraphics[width=0.3\columnwidth]{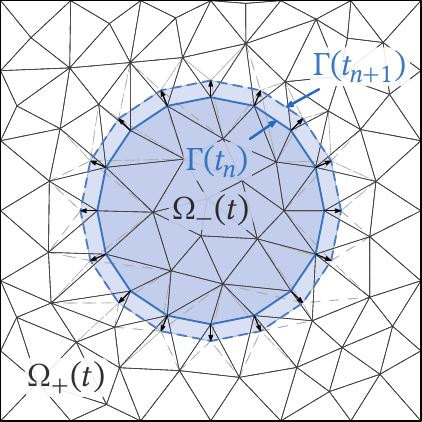}
\caption{Sketch of the two-dimensional moving mesh method, illustrating the interface motion for a time step \(t_n \to t_{n+1}\).}
\label{fig:2d_mmesh}
\end{figure}

The focus of this section lies on the numerical discretization of two-phase flow models of the form \eqref{eq:generic_cl_2p} with a sharp phase boundary. 
This method is utilized as the continuum-scale discretization for our multiscale model. 
For this purpose, we employ a finite volume method (\ipfv{}-method)
that builds upon an 
interface-preserving moving mesh (\ipmm{}) introduced in
\cite{magiera:molecular:2021,alkaemper.magiera.ea:interface:2021}.
In this contribution, we will use the \ipmm{} only as a black box without going into details. 
In summary, the \ipmm{} provides a simplicial conforming moving mesh that 
is able to preserve a discretized interface \(\continterface_{\cT}(t)\).
This means that the interface always coincides with cell surfaces.  
Furthermore, it is possible to prescribe the interface motion, while employing local remeshing (which is described in \cite{alkaemper.magiera.ea:interface:2021}). 
In that way, large mesh deformations are possible, in contrast to e.g. basic Lagrangian approaches. 
The main benefit of the \ipfv{}-method is, that there is no interpolation across the interface, and that the interface dynamics can be directly described by an interface solver.
A graphical illustration of the two-dimensional \ipmm{} is shown in Figure~\ref{fig:2d_mmesh}.
In the following, assume we have a conforming, time-dependent, \ipmm{} \(\cT(t)\), \(t \geq 0\), of the two-phase domain \( \contdomain(t) = \contdomain_\phasel(t) \cup \contdomain_\phaser(t) \cup \continterface(t)  \subset \bR^d\).
Note that the evolution of the interface --- and consequently the moving mesh \(\cT(t)\) --- is not known a priori, but is a result from the interface dynamics during the simulation. 
\newline
We introduce some further notation. 
The set of vertices incident to the interface \(\continterface_{\cT}(t)\) is denoted by \(\cV_{\continterface}(t)\), and the set of interface surfaces by \(\cS_{\continterface}(t)\). 
The cell surfaces that lie inside the bulk phases are given by \(\cS_\phaserl(t) = \cS(t) \mathbin{\backslash} \cS_{\continterface}(t)\), where \(\cS(t)\) denotes the set of all cell surfaces. 
The index set for all vertices  \(\vp \in \cV(t)\) is denoted by \(\bI_{\cV}(t) \subset \bN\), 
and the one for all cells in the mesh \(\cT(t)\) by \(\bI_{\cC}(t) \subset \bN\).

\paragraph*{Numerical Flux}
The \ipfv{}-method is based on a finite-volume scheme. 
As such a numerical flux function
\begin{align}
  G \colon \cU \times \cU \times \bS^{d-1} \to \bR^{m}
\end{align}
is needed. 
For the sake of simplicity, we employ only the Lax--Friedrichs flux function 
\begin{align}
    G(\UU, \vV, \vn) \coloneqq \frac{\vF(\UU; \vn) + \vF(\vV; \vn)}{2} - \alpha_{\mathrm{LF}} (\UU - \vV),
\end{align}
with the parameter \(\alpha_{\mathrm{LF}} > 0\), but it is possible to apply other numerical flux functions that can deal with EOS like \cite{thol.rutkai.ea:equation:2016}.

\paragraph*{The \ipfv{}-Method}
We start by giving a rough outline of the derivation of the \ipfv{}-method. 
For this purpose, we integrate \eqref{eq:generic_cl_2p} over a cell \(C(t) \in \cT(t)\) with time dependent vertices, and the time interval \(t \in [t_n, t_{n+1})\). 
Assuming sufficient regularity, this results after applying Gauss's theorem in
\begin{align}
  \MoveEqLeft
  \int^{t_{n+1}}_{t_n}  \int_{C(t)}  \partial_t \UU \,\ddd \vx \ddd t
  + \int^{t_{n+1}}_{t_n}  \int_{\partial C(t)}  \vf(\UU) \cdot \vn \,\ddd S \ddd t = \vec{0}.
\end{align}
By applying a version of Reynolds' theorem, we arrive at the cell-wise conservative form of  \eqref{eq:generic_cl_2p}, i.e. 
\begin{align}
  \label{eq:ipfv:cellconservative}
  \begin{aligned}    
\int_{C(t_{n+1})}  \UU(\vx, t_{n+1}) \, \ddd \vx 
  ={}& 
  \int_{C(t_{n})}  \UU(\vx, t_{n}) \, \ddd \vx  
- \int^{t_{n+1}}_{t_n}  \int_{\partial C(t)}  \bigl(\vf(\UU) - \UU \vs^\top \bigr) \cdot \vn \, \ddd S \ddd t,
\end{aligned}
\end{align}
where \(\vs \colon \partial C(t) \to \bR^d\) denotes the speed of the surface at a point \(\vec{\xi} \in \partial C(t)\). 
For a thorough derivation we refer to \cite{chalons.rohde.ea:finite:2017}.
\newline
The \ipfv{}-method is basically a discretization of the terms in \eqref{eq:ipfv:cellconservative}, while distinguishing between cell surfaces that form the interface \(\continterface_{\cT}(t)\), and the bulk phase surfaces. 
The algorithmic description follows.

\begin{boxalgorithm}[label=alg:mmfv]{Finite Volume Method on Interface-Preserving Moving Meshes}
\textbf{Input:}
initial mesh \(\cT(0)\), 
cell-averaged initial data \((\UU^0_i)_{i \in \bI_{\cC}(0)}\), 
interface solver \(\genericsolversym \colon \PP_\phasel \times \PP_\phaser \times \bS^{d-1} \to \PP_\phasel \times \PP_\phaser \times \bR\).

\textbf{Parameters:} 
time step \(\Delta t_n > 0\),
vertex motion regularization parameter \(\lambda_{\mathrm{motion}} \geq 0\),
numerical flux function \(G \colon \cU \times \cU \times \bS^{d-1} \to \bR^m\). 

\tcbsubtitle{Algorithm}
For each time step \(t_n \to t_{n+1}\), do
\begin{itemize} 
\item Initialize \(\Delta \UU_i = 0\) for all \(i \in \bI_{\cC}(t_n)\).\item For all interface surfaces \(S^n_{ij} \in \cS_{\continterface}(t_n)\) --- with their incident cells being \(C_i\), \(C_j\) and assuming without loss of generality  \(\UU^n_i \eqqcolon \UU_\phasel \in \cU_\phasel\), \(\UU^n_j \eqqcolon \UU_\phaser \in \cU_\phaser\) --- do 
\begin{itemize} 
 \item Solve the interface Riemann problem 
 \begin{align} \label{eq:mmfv:riemannsolver}
  (\UU^{*}_{\phasel}, \UU^{*}_{\phaser}, s) = \genericsolversym(\UU_\phasel, \UU_\phaser, \vn_{ij}).
 \end{align}
 \item Add contributions to the finite volume update terms
 \begin{align*}
\Delta\UU_i &\minuseq \frac{\Delta t_n \, \abs{S^n_{ij}}}{\abs{C_i(t_n)}} \bigl( f(\UU^{*}_{\phasel}) - s \UU^{*}_{\phasel} \bigr),  
&
\Delta\UU_j &\pluseq \frac{\Delta t_n \, \abs{S^n_{ij}}}{\abs{C_j(t_n)}} \bigl( f(\UU^{*}_{\phaser}) - s \UU^{*}_{\phaser} \bigr).
 \end{align*}
 \item Save \((s, \vn_{ij}, \abs{S_{ij}})\) for every vertex \(\vp_k\) incident to \(S_{ij}\). 
\end{itemize}
\item Update vertex motion \(\vm_k\) for all \(k \in \bI_{\cV}(t_n)\), see \eqref{eq:motion:linear_system}.
\item For every surface inside the bulk domain \(S^n_{ij} \in \cS_\phaserl(t_n)\)
\begin{itemize}
 \item Compute the linear motion contribution
 \begin{align} 
   \label{eq:linear_motion_contribution}
  \ell_{ij} \coloneqq 
  \frac{\abs{S^{n+1}_{ij}} + \abs{S^{n}_{ij}}}{2} 
  \cdot 
   \frac{\UU^{n}_{i} + \UU^{n}_{j}}{2} 
  \cdot 
  \frac{\vn^{n+1}_{ij} + \vn^{n}_{ij}}{2}   
  \cdot 
  \frac{\sum_{ \vp \in \cV(S^n_{ij}) } \vm_k}{\abs{ \cV(S^n_{ij}) }},
 \end{align}
 with \(\cV(S^n_{ij})\) denoting the set of vertices \(\vp \in \cV(t_n)\) that form the cell surface \(S^n_{ij}\).
 Consequently, \(\ell_{ij} = 0\) if \(S_{ij}\) is not incident to an interface vertex.
\item Add contributions to the finite volume update terms
 \begin{align*}
\Delta\UU_i &\minuseq \frac{\Delta t_n}{\abs{C_i(t_n)}} 
\left(
\abs{S^n_{ij}} \, G(\UU^n_i, \UU^n_j, \vn^n_{ij}) 
+ \ell_{ij} \right),
\\
\Delta\UU_j &\pluseq \frac{\Delta t_n}{\abs{C_j(t_n)}} 
\left(
\abs{S^n_{ij}} \, G(\UU^n_i, \UU^n_j, \vn^n_{ij}) 
+ \ell_{ij} \right)
.
 \end{align*}
\end{itemize}
\item Add update terms 
to all cells: 
\begin{align*}
\UU^{n+1}_i = \UU^{n}_i + \Delta \UU_i,
\text{ for } i \in \bI_{\cC}(t_n).
\end{align*}
\item Update vertex positions according to 
\begin{align}
  \label{eq:ipfvm:vertex_motion}
\vp^{n+1}_k = \vp^{n}_k + \Delta t_n \, \vm_k, 
\text{ for } k \in \bI_{\cV}(t_n),
\end{align}
and rescale the data of all neighboring cells according to the changing volume, i.e.
\begin{align}
  \label{eq:motion_data_scaling}
  \UU_i \mapsfrom  \frac{\abs{C_i(t_n)}}{\abs{C_i(t_{n+1})}}  \UU_i.
\end{align}
\item Perform the remeshing according to the \ipmm{} (see \cite{alkaemper.magiera.ea:interface:2021}).
\end{itemize}
\end{boxalgorithm}
\begin{remark}
Concerning the \ipfv{}-method in Algorithm~\ref{alg:mmfv}, there are some remarks to be made.
\begin{itemize}
  \item 
  The finite volume update terms \(\Delta \UU_i\) are used as variables, with \(\pluseq\), \(\minuseq\) denoting assignment by addition or subtraction.
  \item 
  By scaling the data in \eqref{eq:motion_data_scaling} the interface motion remains conservative. 
  Even more so, we avoid mixing phase states in the context of two-phase flows. 
  \item 
  The linear motion contribution \eqref{eq:linear_motion_contribution} originates from discretizing the motion-dependent part of the flux term in \eqref{eq:ipfv:cellconservative} for cell surfaces adjacent to the interface. 
\end{itemize}
\end{remark}

It remains to discuss how the interface motion is computed. 

\paragraph*{Computing the Vertex Motion of the Mesh}
In the \ipfv{}-scheme (Algorithm~\ref{alg:mmfv}) 
the interface speed \(s\) is always computed for each cell surface \(S \in \cS_{\continterface}(t)\) in direction of its normal vector \(\vn\) --- however, for the \ipmm{} we require the motion at each vertex. 
As such we have to compute the vertex motion \(\vm_k = \vm_k(t) \in \bR^d\) for each interface vertex \(\vp_k \in \cV_{\continterface}(t)\). 
Let us assume that for each interface surface \(S_l \in \cS_{\continterface}(t)\) incident to \(\vp_k\) we have already computed the interface wave speed \(s_l\) in normal direction \(\vn_l\).   
In the multiscale model, this will be provided by the interface solver.

A natural choice for the interface motion \(\vm_k\) of \(\vp_k\) is one that satisfies
\begin{align}
\label{eq:motion_projection}
 \vn_l \cdot \vm_k = s_l \text{ for all surfaces } S_l \text{ incident to } \vp_k, 
\end{align}
which means that the projection of \(\vm_k\) in each normal direction \(\vn_l\) coincides with the interface speed \(s_l\).
\newline
In one space dimension this is trivially fulfilled by setting \(\vm_k = s\). 
In case of two space dimension, one vertex has usually (in our setting) only two interface edges, i.e. \eqref{eq:motion_projection} boils down to a \(2 \times 2\) system of linear equations, that is solvable if the interface edges are not parallel. 
In three space dimensions each interface vertex has in most cases more than three incident interface surfaces, consequently, the linear system \eqref{eq:motion_projection} is generally overdetermined.  
To solve this problem, linear least squares fitting, using the normal equations, can be performed.
For that, we define 
\begin{align*}
 \mat{N} \coloneqq \begin{pmatrix}
                    \vn^\top_1 \\ \vdots \\ \vn^\top_K
                   \end{pmatrix}, 
 \quad 
 \vs \coloneqq \begin{pmatrix}
                    s_1 \\ \vdots \\ s_K
               \end{pmatrix}, 
 \end{align*}
where \(\vn^\top_l\), \(s_l\) correspond to the normal vectors and surface speeds for all \(K \in \bN\) incident surfaces \(S_l\) of \(\vp_k\). 
Then the normal equations are given by 
\begin{align}
  \mat{N}^\top \mat{N} \vm_k = \mat{N}^\top \vs.
\end{align}
To regularize this system, we apply the generalized Tikhonov regularization \cite{tichonov.leonov.ea:nonlinear:1998}, with the expected value \(\overline{\vm}_k\) being the averaged motion 
\begin{align} 
\label{eq:avg_motion}
\overline{\vm}_k \coloneqq \left(\dfrac{1}{\abs{\cS_{\continterface}(\vp_k)}}
\sum_{S_l  \in \cS_{\continterface}(\vp_k) } s_l\right) \dfrac{\overline{\vn}_k}{\norm{\overline{\vn}_k}}, 
\end{align}
with \(\cS_{\continterface}(\vp_k)\) denoting the interface surfaces incident to the vertex \(\vp_k\), and the averaged normal vector given by 
\begin{align}
\overline{\vn}_k \coloneqq \dfrac{1}{\abs{\cS_{\continterface}(\vp_k)}}
\sum_{S_l  \in \cS_{\continterface}(\vp_k) } \vn_l.
\end{align}
Then, the motion \(\vm_k\) of the vertex \(\vp_k\) is the solution of the linear system 
\begin{align}
  \label{eq:motion:linear_system}
  (\mat{N}^\top \mat{N} + \lambda_{\mathrm{motion}} \mat{I}) \vm_k = \mat{N}^\top \vs + \lambda_{\mathrm{motion}} \overline{\vm}_k,
\end{align}
where \(\lambda_{\mathrm{motion}} > 0\) is the regularization parameter, and \(\mat{I}\) the \((d\times d)\)-identity matrix.

\begin{remark}
It is noteworthy that the \ipfv{}-method is applicable beyond the multiscale model presented in this work. 
Many models can be applied in the bulk domains or the interface, and the method is not restricted to phase boundary scenarios.
An open-source implementation of the \ipmm{} can be found at \cite{alkaemper:interface:2022}. 
\end{remark}
 
\section{The Multiscale Model}
\label{sec:multiscale}

In this core section we present the multiscale model and the interplay of each of its components.
We first describe it on the general level of Section~\ref{sec:generic_two_phase}. 
Thereafter, we will go into details and explain the specific implementation for liquid--vapor flow (Section~\ref{sec:euler}), where the interface dynamics are described by \mdshort{} simulations. 

\subsection{The Multiscale Model for General Two-Phase Flow}

\begin{figure}[tbp]
\centering
\includegraphics[width=0.75\columnwidth]{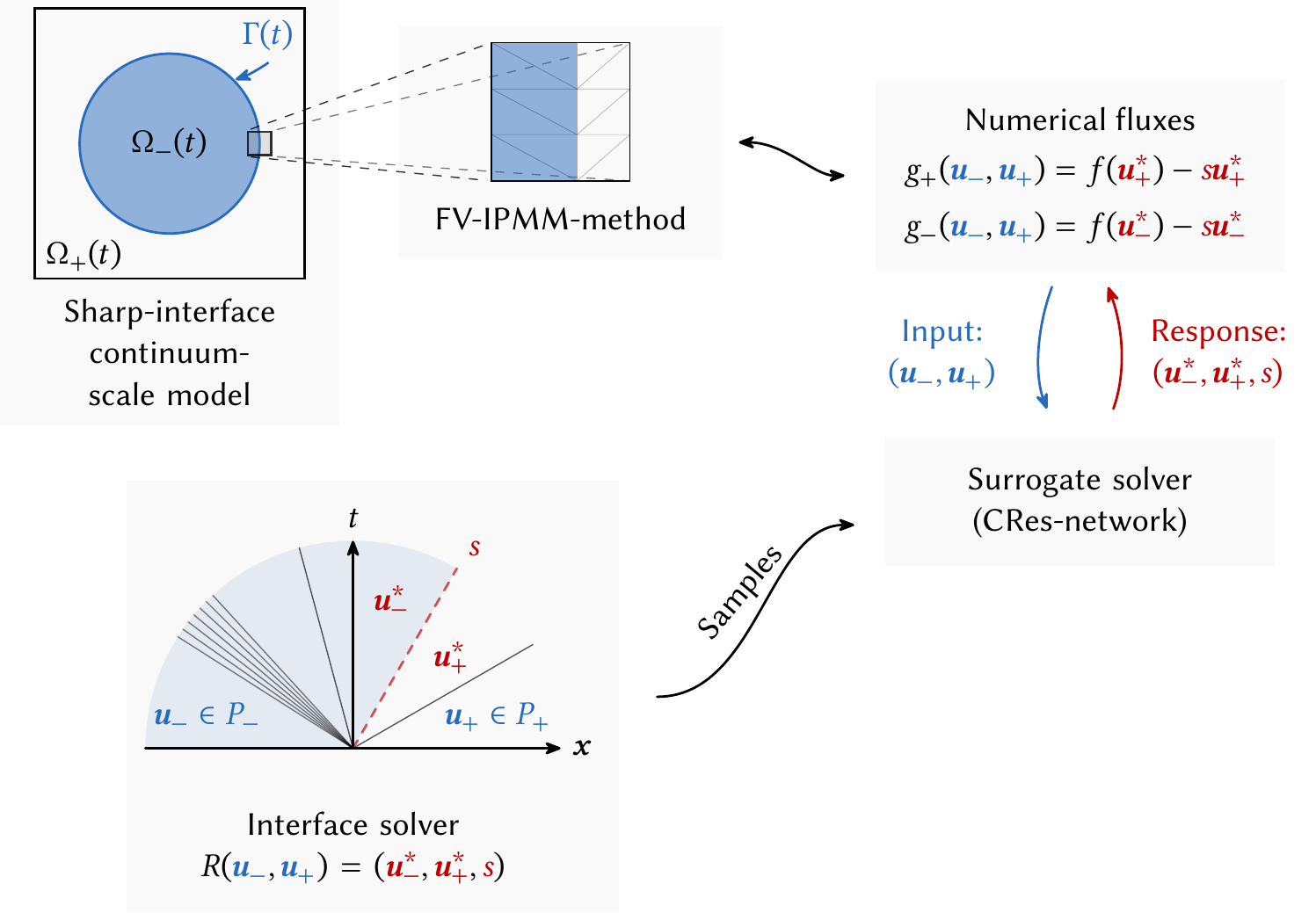}
\caption{
Sketch of the multiscale model, its main components, and their interactions. 
}
\label{fig:multiscale_scheme}
\end{figure}

The multiscale model consists of several components.
On the \emph{continuum scale}, two-phase fluid flow can be modeled by conservation laws with sharp phase boundaries, see Section~\ref{sec:generic_cl}. 
Models of this type can be discretized via the 
\emph{\ipfv{}-method} described in Section~\ref{sec:ipvm_algorithm}. 
A major benefit of this approach is, that it permits us to model the fluid dynamics directly at the moving interface, for example, by prescribing a general interface solver \(\genericsolversym\) (see Definition~\ref{def:generic_interface_solver}).
Following this idea, we model the interface dynamics with a multiscale model that is based on microscale simulations of the phase boundary, forming the interface solver. 
Such a \emph{microscale interface solver} \(\rmicro\) has the form given in  \eqref{eq:riemannsolver:rotated}.
For given initial states 
\(\uu_{\phasel} \in \pp_{\phasel}\), \(\uu_{\phaser} \in \pp_{\phaser}\)
it solves the microscale Riemann problem (in normal direction) and extracts the wave speed \(s \in \bR\) of the interface as well as the adjacent wave states 
\(\uu^*_{\phasel} \in \pp_{\phasel}\), \(\uu^*_{\phaser} \in \pp_{\phaser}\) --- see Figure~\ref{fig:generic_wave_pattern} for an exemplary sketch of corresponding wave pattern. 
The microscale interface solver \(\rmicro\) can be used to formulate a general interface solver \(\genericsolversym\)  that is used in the \ipfv{}-method. 
\newline
Specifically, for liquid--vapor flow, the microscale interface solver \(\rmicro\) consists of \mdshort{} simulations of  molecular-scale Riemann problems. 
From these particle simulations the phase boundary dynamics are inferred. 
The \mdshort{} interface solver is described in Section~\ref{sec:md}.
\newline
A graphical representation of the interplay between the different components of the multiscale model is shown in Figure~\ref{fig:multiscale_scheme}.

\begin{remark}
  \label{remark:different_interface_solvers}
  The multiscale model is not restricted to employ particle-based microscale solvers for the general interface solver \(\genericsolversym\). 
  There are many other possibilities that can be used. 
  For example one could employ exact Riemann solvers \cite{chen.hattori:exact:2015,zeiler:liquid:2015,schleper:hll:2016}, or (computationally cheaper) relaxation solvers \cite{fechter.jaegle.ea:exact:2013}. 
  The interface solver may also be based on other microscale models like diffuse-interface approximations \cite{blesgen:generalization:1999} or completely different modeling approaches.
\end{remark}

\subsubsection{The Multiscale Model Algorithm}

In this section, the application of the multiscale model is explained.
It employs the surrogate interface solver \(\networkF_{\nnparams}\) in the \ipfv{}-discretization (Section~\ref{sec:ipvm_algorithm}) of a continuum-scale two-phase model. 

We start with the initial data on the continuum scale, which means, we prescribe the initial bulk phase domains \(\contdomain_{\phaserl}(0)\), the
initial interface position \(\continterface(0)\), and the corresponding initial conditions \(\UU_{\phaserl,0}(\vx)\). 
Both domains \(\contdomain_{\phaserl}(0)\) are triangulated while approximating the interface \(\continterface(0)\).
Then, we employ the \ipfv{}-method (see Algorithm~\ref{alg:mmfv}).
In the following, the multiscale model is given in algorithmic form.

\begin{boxalgorithm}[label=alg:multiscale]{Multiscale Model Algorithm}
\textbf{Input:} 
initial bulk domains \(\contdomain_{\phaserl}(0)\), 
initial interface position \(\continterface(0)\), 
and initial conditions \(\UU^{\phaserl}_{0}(\vx)\).

\textbf{Prerequisite:} 
a surrogate solver \(\networkF_{\nnparams}\), that is based on the microscale interface solver \(\cR\). 
\tcbsubtitle{Algorithm}
\begin{itemize} 
  \item Approximate \(\contdomain_{\phaserl}(0)\) and interface \(\continterface(0)\) by mesh with interface \(\continterface_{\cT}\), 
  comprising the initial moving mesh \(\cT(0)\).
\item Compute cell-averaged data \((\UU^0_i)_{i \in \bI_{\cC}(0)}\) on \(\cT(0)\), using the initial conditions \(\UU_{\phaserl,0}(\vx)\).
\item Set up the interface solver 
  \begin{align}
    \label{eq:multiscale:algorithm:interface_solver}
    \genericsolversym \colon \PP_{\phasel} \times \PP_{\phaser} \times \bS^{d-1} \to \PP_{\phasel} \times \PP_{\phaser} \times \bR
    :
    (\UU_{\phasel}, \UU_{\phaser}, \vn) \mapsto (\UU^{*}_{\phasel}, \UU^{*}_{\phaser}, s),    
  \end{align}
  where the mapping \( (\UU_{\phasel}, \UU_{\phaser}, \vn) \mapsto (\UU^{*}_{\phasel}, \UU^{*}_{\phaser}, s)\) is defined by the following steps:
  \begin{enumerate}
    \item Map the states \(\UU_{\phaserl}\) in direction of \(\vn\), according to \eqref{eq:directional_states}, resulting in \(\uu_{\phaserl} = \UU_{\parallel\vn, \phaserl}\) and \(\UU_{\perp\vn, \phaserl}\).    
     \item 
     \label{item:multiscale:apply_surrogate_solver}
     Evaluate the surrogate solver (for \(\rmicro\))
     \begin{align}
       (\uu^{*}_{\parallel\vn, \phasel}, \uu^{*}_{\parallel\vn, \phaser}, s) = \networkF_{\nnparams}(\uu_{\phasel}, \uu_{\phaser}).
     \end{align}    
     \item Project the directional states: 
     \begin{align}
       \UU^{*}_{\phaserl}  
      = P_{\vn}(\uu^{*}_{\phaserl})  + \UU_{\perp\vn, \phaserl},
     \end{align}  
     with \(P_{\vn}\) defined as in \eqref{eq:reverse_projection}.
\end{enumerate}  
  \item Run the \ipfv{}-method (Algorithm~\ref{alg:mmfv}) with the initial mesh \(\cT(0)\), the 
   initial data \((\UU^0_i)_{i \in \bI_{\cC}(0)}\), and using the 
  interface solver \(\genericsolversym\) defined in \eqref{eq:multiscale:algorithm:interface_solver}.  
\end{itemize}
\textbf{Result:} 
approximate solution \(\UU_{\cT}(\vx,t) = U^{n}_{i}\), for \(\vx \in C_i(t_n)\) and \(t \in [t_n, t_{n+1})\), on the moving mesh \(\cT(t)\).
\end{boxalgorithm}

\begin{remark}[Active Sampling Strategies]  
\label{remark:active_sampling}
One major drawback of a fixed surrogate solver is, that it is restricted to a predetermined flow regime, which is due to the static bounding region \(B_{\mathrm{in}}\).
To mitigate this problem, it is possible to employ an active sampling strategy, which in this context is presented in our work \cite{magiera.rohde:particle:2018}. 
The main idea is to introduce a score value for an arbitrary input \(\datx\), that measures the uncertainty of the surrogate solver on this specific input point. 
If the uncertainty becomes too large, the original solver is evaluated, and the resulting sample \( \daty = \cR(\datx) \) is added to the training data set and the surrogate solver gets retrained. 
\end{remark}

 \subsection{The Multiscale Model for Liquid--Vapor Flow}
\label{sec:multiscale_nonisothermal}
In this section, we apply the multiscale model to compressible liquid--vapor flow.
On the continuum scale, the two-phase flow is described by the Euler equations, as in Section~\ref{sec:euler}. 
The interface dynamics on the microscale are modeled by \mdshort{} simulations, see Section~\ref{sec:md}, which are substituted by a surrogate solver, see Section~\ref{sec:surrogate}.

\paragraph*{Microscale Interface Solver}
To define the liquid--vapor flow  multiscale model, we have to specify a microscale interface solver of the form
\begin{align}
  \label{eq:multiscale_noniso_riemannsolver}
  \rmicro \colon 
  (\rho_\phasel, \vm_\phasel, E_\phasel,  \rho_\phaser, \vm_\phaser, E_\phaser, \vn) 
  \mapsto 
  (\rho^*_\phasel, \vm^*_\phasel, E^*_\phasel, \rho^*_\phaser, \vm^*_\phaser, E^*_\phaser, s).
\end{align}
The solver \eqref{eq:multiscale_noniso_riemannsolver} is based on \mdshort{} simulations (Section~\ref{sec:md})
and is realized as follows.
Let the  initial states \((\rho_\phaser, \vm_\phaser, E_\phaser) \in \PP_{\phaser}\), \((\rho_\phasel, \vm_\phasel, E_\phasel) \in \PP_{\phasel}\), and normal direction \(\vn \in \bS^{d-1}\) be given. 
We compute the specific inner energies \(\varepsilon_{\phaserl} = \tfrac{E_{\phaserl}}{\rho_{\phaserl}} - \tfrac{\norm{\vm_{\phaserl}}^2}{2 \rho^2_{\phaserl}}\), and temperatures \(T_{\phaserl} = T(\rho_\phaserl, \varepsilon_{\phaserl})\), using the \eosshort{}.
Next, map \((\rho_\phaserl, \vm_\phaserl, T_\phaserl)\) in direction of \(\vn\) according to 
\eqref{eq:directional_states}, i.e. 
\begin{align}
  \label{eq:noniso:multiscale:projected_state_variables}
  \uu_{\phaserl} =
  (\rho_\phaserl, m_\phaserl, T_\phaserl)
  &\coloneqq 
  \UU_{\parallel\vn, \phaserl} 
  = (\rho_\phaserl, \vm_\phaserl \cdot \vn, T_\phaserl), 
  \\
\UU_{\perp\vn, \phaserl} 
  &\coloneqq 
  (0, \vm_\phaserl - (\vm_\phaserl \cdot \vn) \vn, 0 ).
\end{align}    
Then, we consider the fluid states in the  reference frame of Section~\ref{sec:surrogate:liquid-vapor}, which means we choose the \mbox{\(\phasel\)-phase} velocity \(v_{\phasel} = \frac{m_\phasel}{\rho_\phasel}\) as the reference velocity \(v_{\mathrm{ref}} \coloneqq v_\phasel \).
Thus, we consider the fluid velocities
\(\widetilde{v}_\phaserl = v_{\phaserl}  - v_{\mathrm{ref}}\).
Next, we  solve the Riemann problem \eqref{eq:rotated_cl} in normal direction \(\vn\) with Riemann data \((\rho_\phaserl, \rho_\phaserl\widetilde{v}_\phaserl, T_{\phaserl})\) using the microscale model.
This means, we run a \mdshort{} simulation (Algorithm~\ref{alg:md_riemann}) for the initial continuum states \((\rho_\phasel, \widetilde{v}_\phasel, T_{\phasel})\) and \((\rho_\phaser, \widetilde{v}_\phaser, T_{\phaser})\),
i.e. evaluating \(\rmd\) from \eqref{eq:md_noniso_riemann_solver}
resulting in the  states \((\rho^*_\phasel, \widetilde{m}^*_\phasel, T^*_\phasel)\), \((\rho^*_\phaser, \widetilde{m}^*_\phaser, T^*_\phaser)\), and the wave speed \(\widetilde{s}\) in normal direction \(\vn\), 
with \(\widetilde{m}^*_\phaserl = \rho^*_\phaserl \widetilde{v}^*_\phaserl\).
Next, we return the reference velocity, i.e. 
\begin{align}
  m^*_\phaserl = \widetilde{m}^*_\phaserl + \rho^*_{\phaserl} v_{\mathrm{ref}}, 
  \qquad
  s = \widetilde{s}  + v_{\mathrm{ref}}.
\end{align} 
The directional states are projected back to the full-dimensional states
\begin{align}
  (\rho^*_\phaserl, \vm^*_\phaserl, T^*_\phaserl)     
 = (\rho^*_\phaserl, m^*_\phaserl \vn, T^*_\phaserl)  + \UU_{\perp\vn, \phaserl}
 = P_{\vn}((\rho^*_\phaserl, m^*_\phaserl, T^*_\phaserl))  + \UU_{\perp\vn, \phaserl},
\end{align}  
with \(P_{\vn}\) defined as in \eqref{eq:reverse_projection}.
After computing the specific inner energies \(\varepsilon^*_{\phaserl} = \varepsilon(\rho^*_{\phaserl}, T^*_{\phaserl})\), and energy densities \(E^*_{\phaserl} = \rho^*_{\phaserl} \varepsilon^*_{\phaserl} + \tfrac{\norm{\vm^*_{\phaserl}}^2}{2 \rho^*_{\phaserl}}\), we have obtained the 
wave states \((\rho^*_\phaser, \vm^*_\phaser, E^*_{\phaser})\), \((\rho^*_\phasel, \vm^*_\phasel, E^*_{\phasel})\), and the interface speed \(s\) in direction of \(\vn\).
 
\section{Numerical Simulation Results for the Liquid--Vapor Multiscale Model}
\label{sec:results:noniso}

In this section, we present simulation results for the liquid--vapor multiscale model, as presented in Section~\ref{sec:multiscale_nonisothermal}. 
We consider the two-phase Euler equations of Section~\ref{sec:euler} on the continuum scale, and \mdshort{} simulations for the phase boundary dynamics, as described in Section~\ref{sec:md}.
The parameters used for the simulations are found in 
Appendix~\ref{appendix:parameter_table:noniso}.

\subsection{One-dimensional Results}
\label{sec:results:noniso:1d}

To start with, we consider the multiscale model in one space dimension. 
Our focus lies hereby on Riemann problems, which are compared with their corresponding \mdshort{} simulation result. 
In that way, we are able to validate the consistency 
of the multiscale model with respect to the microscale \mdshort{} simulations.
We consider two settings: 
a pressure-driven wave and a vapor wave colliding with liquid fluid.

\begin{figure}[tp]
\centering
\includegraphics[width=0.75\columnwidth]{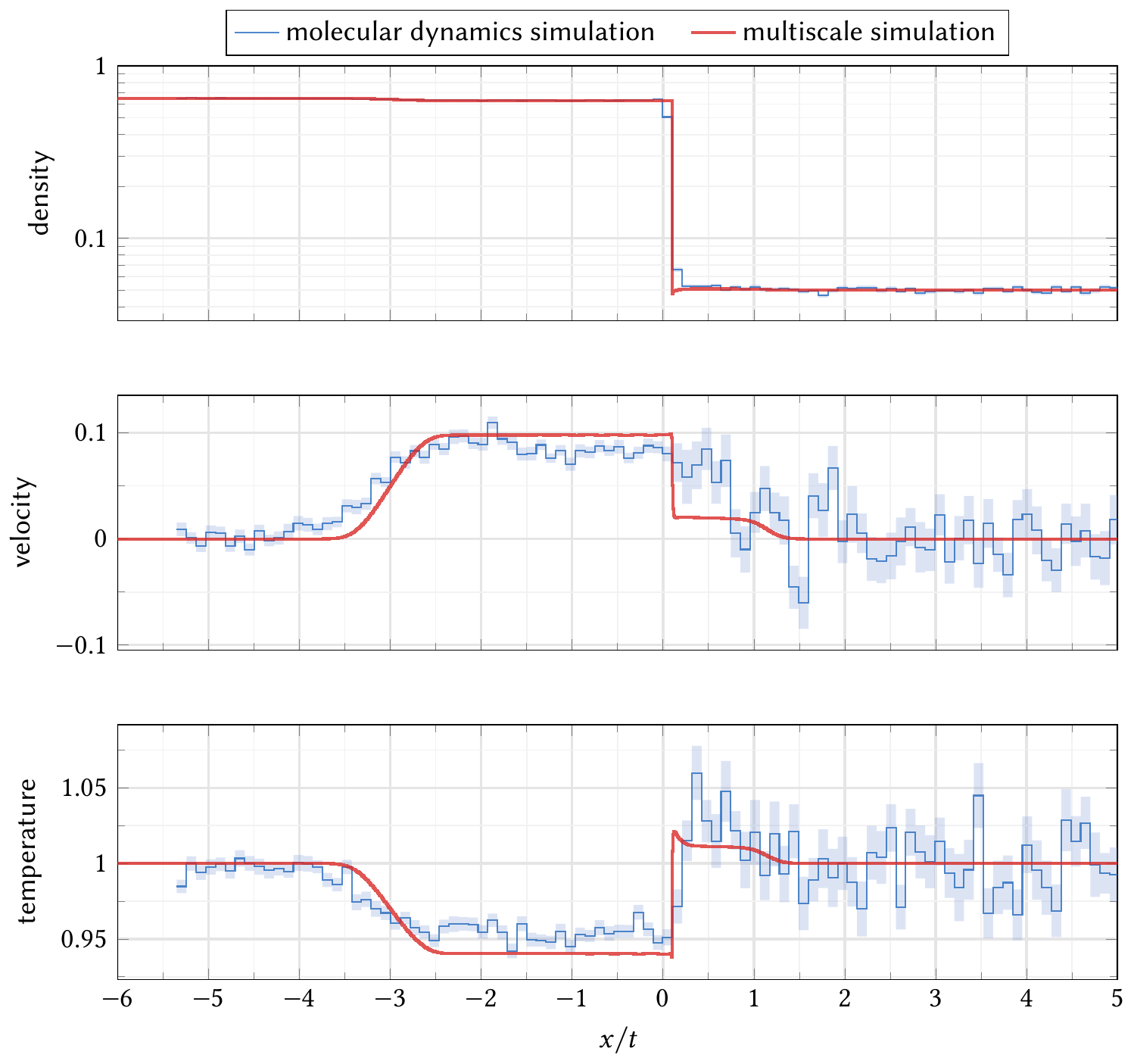}
\caption{ 
Multiscale simulation for the liquid--vapor flow model in one space dimension, overlaid over the corresponding \mdshort{} simulations (averaged over \num{50} simulations and \num{200} bins for local averaging, the shaded region indicates the standard error of the mean). 
The initial conditions 
\(\rho_\phasel = \num{0.65}\), 
\(T_\phasel = \num{1.0}\), 
\(\rho_\phaser = \num{0.05}\), 
\(T_\phaser = \num{1.0}\),
\(v_\phaserl = \num{0}\)
 yield a pressure-driven wave that expands from the liquid phase. 
}
\label{fig:mm1d_noniso_shock}
\end{figure}

\paragraph*{Pressure-Driven Wave}
For the first example, we consider the Riemann initial data 
\(\rho_\phasel = \num{0.65}\), 
\(T_\phasel = \num{1.0}\), 
\(\rho_\phaser = \num{0.05}\), 
\(T_\phaser = \num{1.0}\),
\(v_\phaserl = \num{0}\). 
As seen in Figure~\ref{fig:mm1d_noniso_shock}, this results in a pressure-driven wave, expanding the liquid phase. 
Due to the expansion and the decreasing density inside the liquid phase, the temperature decreases as well. 
The vapor next to the interface is slightly compressed and therefore  has a slightly increased temperature. 
Some mass is transferred from the vapor into the liquid phase, according to the (averaged) relative mass flux \(j \approx \num[round-mode=places,round-precision=3]{-0.004484}\).
\newline
The \mdshort{} simulation results are close to the multiscale solution. 
Specifically, this holds for the speed of the observed waves.
The wave plateau values for the density coincide very well --- however, for the velocity and the temperature, there is some discrepancy. 
This is due to the higher fluctuation in statistical averaging on the vapor side, as on the molecular-scale fewer particles are present there.

\begin{figure}[tp]
\centering
\includegraphics[width=0.75\columnwidth]{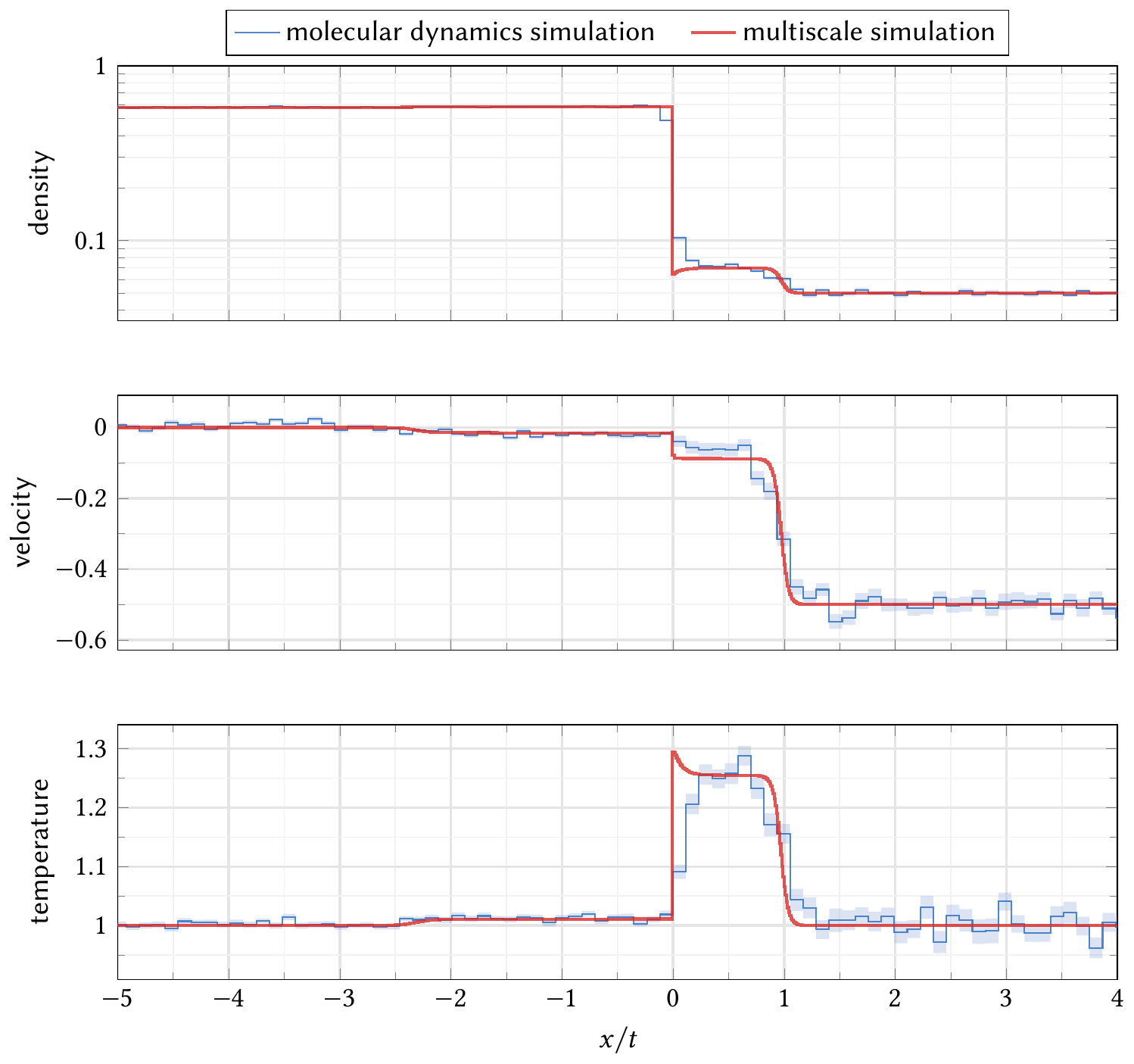}
\caption{ 
Multiscale simulation for the liquid--vapor flow model in one space dimension and the corresponding \mdshort{} simulations (averaged over \num{50} simulations and \num{200} bins for local averaging, the shaded region indicates the standard error of the mean). 
The initial conditions 
\(\rho_\phasel = \num{0.58}\), 
\(v_\phasel = \num{0}\), 
\(T_\phasel = \num{1.0}\), 
\(\rho_\phaser = \num{0.05}\), 
\(v_\phaser = \num{-0.5}\), 
\(T_\phaser = \num{1.0}\)
are chosen in such a way that a vapor phase wave hits the liquid phase. 
}
\label{fig:mm1d_noniso_condensation}
\end{figure}

\paragraph*{Liquid--Vapor Collision}
In the second example, we consider a vapor wave that hits the liquid phase.
This is implemented by the Riemann initial data 
\(\rho_\phasel = \num{0.58}\), 
\(v_\phasel = \num{0}\), 
\(T_\phasel = \num{1.0}\), 
\(\rho_\phaser = \num{0.05}\), 
\(v_\phaser = \num{-0.5}\), 
\(T_\phaser = \num{1.0}\). 
The simulation result is shown in Figure~\ref{fig:mm1d_noniso_condensation}.
In the vapor phase we observe an increase in the density and temperature. 
The averaged relative mass flux \(j \approx \num[round-mode=places,round-precision=3]{-0.006092}\) indicates that mass is transferred into the liquid phase. 
In the same way, the momentum of the vapor wave is transferred into the liquid phase as seen in the velocity plot. 
\newline
Most noteworthy is the excellent correspondence between multiscale solution and \mdshort{} simulation --- only the temperature distribution of the \mdshort{} simulation exhibits some heat conductivity effects which are not present on the continuum scale. 

In summary, the one-dimensional simulation results are consistent with the molecular-scale simulations, which validates the multiscale approach.

\subsection{Two-dimensional Results}
\label{sec:results:noniso:2d}

Next, we show the performance of the multiscale model in two space dimensions.
For this scenario we consider the domain \(\contdomain = (-1.5, 1.5)^2\), with the two initial phase domains \(\contdomain_\phasel(0) = \{ \vx = (x_1, x_2) \in \bR^2 \setsep \norm{\vx}^2_2 < 0.15\}\) and  \(\contdomain_\phaser(0) = \contdomain \mathbin{\backslash} \overline{\contdomain_\phasel(0)}\).  
The initial data in each subdomain is set to
\begin{align}
  (\rho, \vv, T)(\vx, t)
  = 
  \begin{cases}
    (0.62, (0, 0), 0.8), & \vx \in \contdomain_\phasel(0), \\
    (0.06, (0, 0), 0.9), & \vx \in \contdomain_\phaser(0) \text{ and } x_1 \geq  -0.5, \\
    (0.06, (0.5, 0), 0.9), & \vx \in \contdomain_\phaser(0) \text{ and }  x_1 < -0.5.
  \end{cases}
\end{align}
This corresponds to a vapor wave colliding with a liquid droplet. 
At the left boundary, at \(x_1 = -1.5\), we apply inflow boundary conditions by setting ghost cell values to \((\rho, \vv, T) = (0.06,\allowbreak (0.5, 0.0),\allowbreak 0.9)\). 
Everywhere else we use outflow boundary conditions by copying the inner cell values to the ghost cells.

\begin{figure}[tbp]
 \centering
\includegraphics[width=0.7\columnwidth]{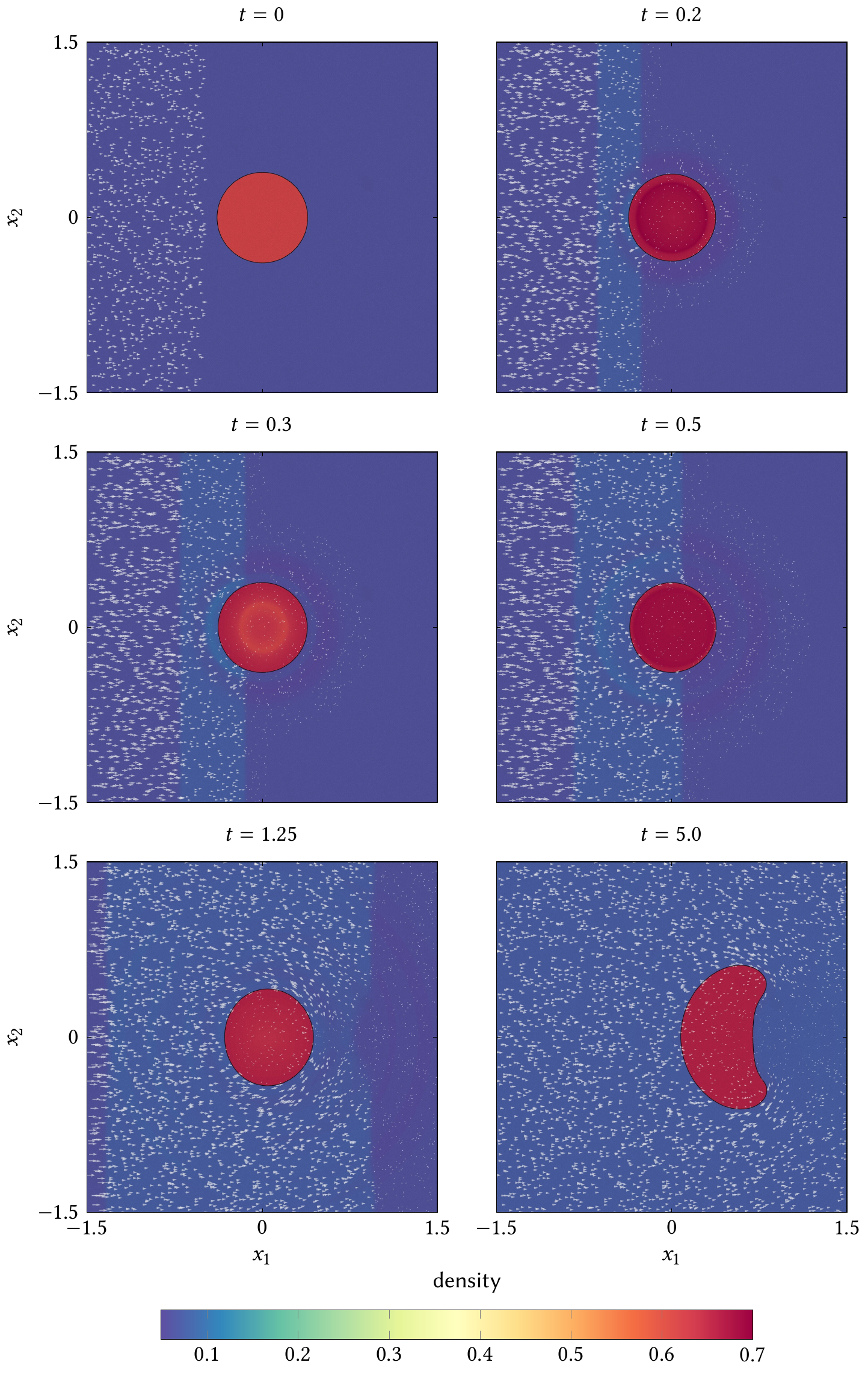}
\caption{
Two-dimensional multiscale simulation of the liquid--vapor flow model. 
The sub-figures depict the density and velocity of the solution at various time steps. 
The corresponding plots of the temperature can be found in Figure~\ref{fig:noniso:shock_droplet:temperature:2d}.
 }
 \label{fig:noniso:shock_droplet:density:2d}
\end{figure}

\begin{figure}[tbp]
\centering
\includegraphics[width=0.7\columnwidth]{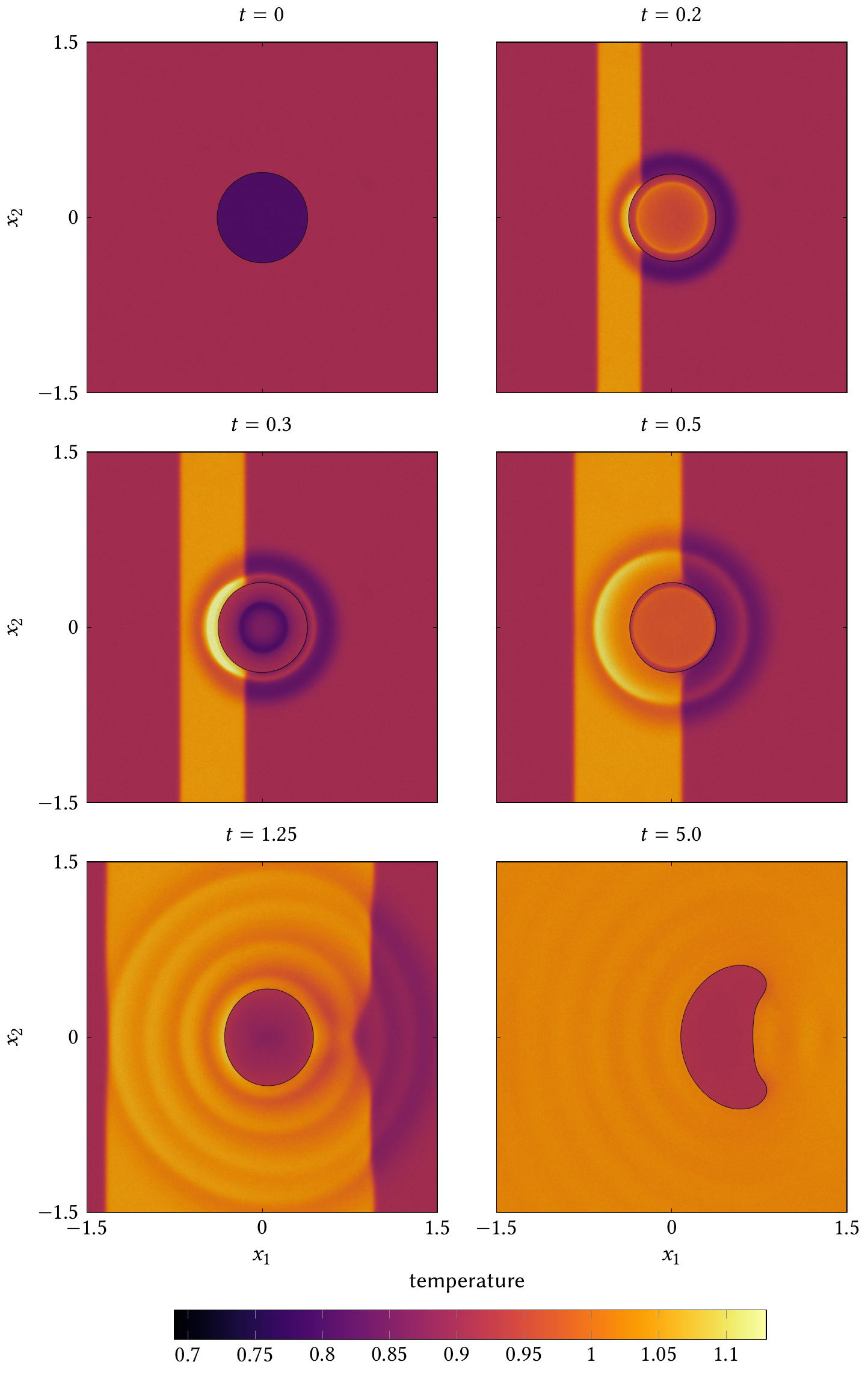}
\caption{
 Two-dimensional multiscale simulation of the liquid--vapor flow model. 
 The sub-figures depict the temperature of the solution at various time steps. 
 The corresponding plots of the density can be found in Figure~\ref{fig:noniso:shock_droplet:temperature:2d}.
}
\label{fig:noniso:shock_droplet:temperature:2d}
\end{figure}

\begin{figure}[tbp]
\centering
\includegraphics[width=0.6\columnwidth]{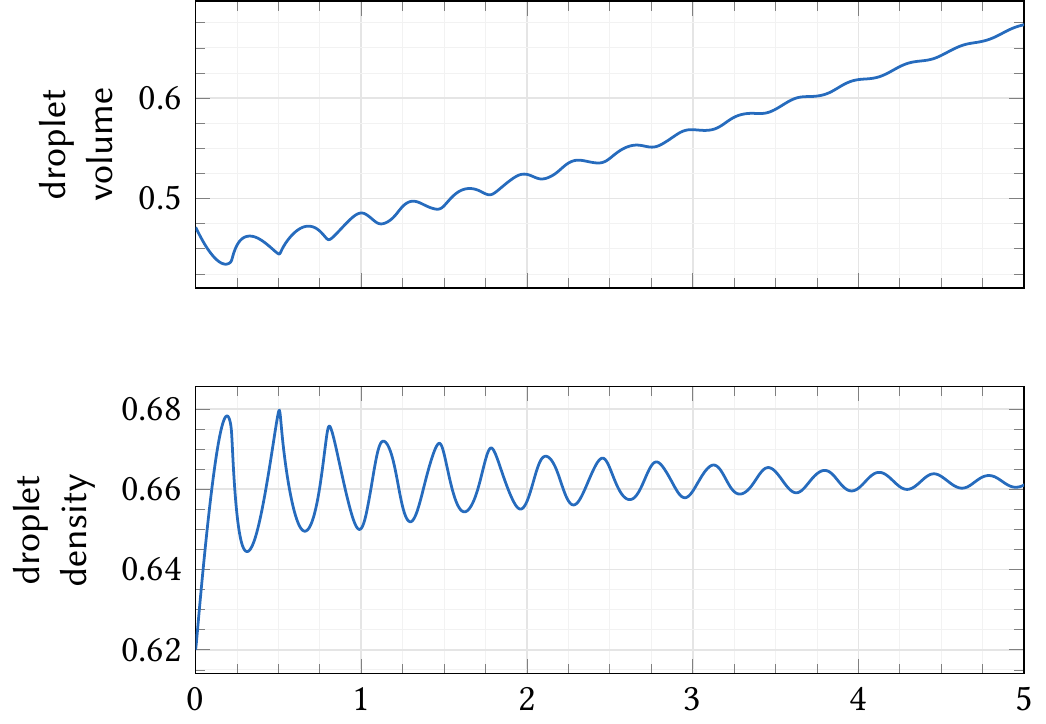}
\caption{ 
Time evolution of the volume of the droplet and the  droplet-averaged density for the liquid--vapor flow multiscale simulation in two space dimensions.
} 
\label{fig:2d_noniso_data_averages}
\end{figure}

The density and velocity field of the simulation are shown in Figure~\ref{fig:noniso:shock_droplet:density:2d} for several time steps; the corresponding temperature field is depicted in Figure~\ref{fig:noniso:shock_droplet:temperature:2d}.
As the initial data does not coincide with equilibrium states, the liquid droplet starts to oscillate --- clearly seen in Figure~\ref{fig:2d_noniso_data_averages}. 
Meanwhile, the vapor wave increases the density and temperature inside the vapor phase. 
Then, the vapor wave interacts with the droplet, deforms it, and transports it through the domain. 
Note that surface tension is not included in the model. 
Therefore, the droplet deformation might be exaggerated.  

\subsection{Three-dimensional Results}
\label{sec:results:noniso:3d}

In this section we demonstrate the capabilities of the multiscale model in three space dimensions. 
For this setting we consider the domain \(\contdomain = (-1.5, 1.5)^3\), with the two initial phase domains \(\contdomain_\phasel(0) = \{ \vx = (x_1, x_2, x_3) \in \bR^3 \setsep \norm{\vx}^2_2 < 0.15\}\) and  \(\contdomain_\phaser(0) = \contdomain \mathbin{\backslash} \overline{\contdomain_\phasel(0)}\).  
For the initial data in each subdomain we set 
\begin{align}
  (\rho, \vv, T)(\vx, t)
  = 
  \begin{cases}
    (0.62, (0, 0, 0), 0.8), & \vx \in \contdomain_\phasel(0), \\
    (0.06, (0, 0, 0), 0.9), & \vx \in \contdomain_\phaser(0) \text{ and } x_1 \geq  -0.5, \\
    (0.06, (0.5, 0, 0), 0.9), & \vx \in \contdomain_\phaser(0) \text{ and } x_1 < -0.5.
  \end{cases}
\end{align}
For the boundary conditions we apply inflow boundary conditions at 
\(x_1 = -0.5\) with the ghost cell value \((\rho, \vv, T) = (0.06, \allowbreak (0.5, 0, 0), \allowbreak 0.9)\). 
Everywhere else outflow boundary conditions are applied.
This corresponds to the two-dimensional setting of Section~\ref{sec:results:noniso:2d}.

\begin{figure}[tbp]
 \centering
\includegraphics[width=0.7\columnwidth]{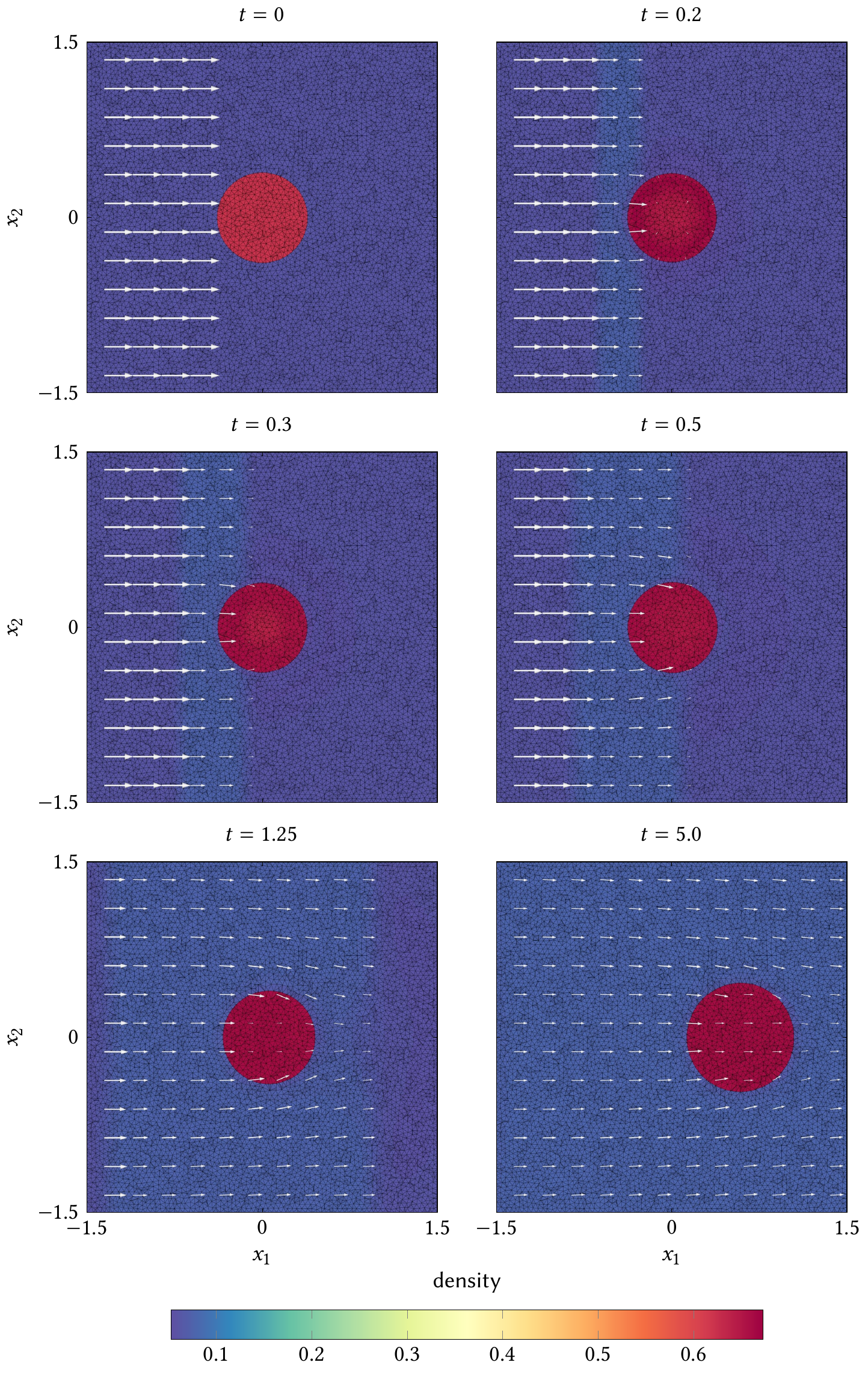}
\caption{
 Slice through the three-dimensional multiscale simulation of the liquid--vapor flow model at \(x_3 = 0\). 
 The sub-figures depict the density and velocity of the solution at various time steps. 
 The corresponding plots of the temperature can be found in Figure~\ref{fig:noniso:shock_droplet:temperature:3d}.
  }
 \label{fig:noniso:shock_droplet:density:3d}
\end{figure}

\begin{figure}[tbp]
 \centering
\includegraphics[width=0.7\columnwidth]{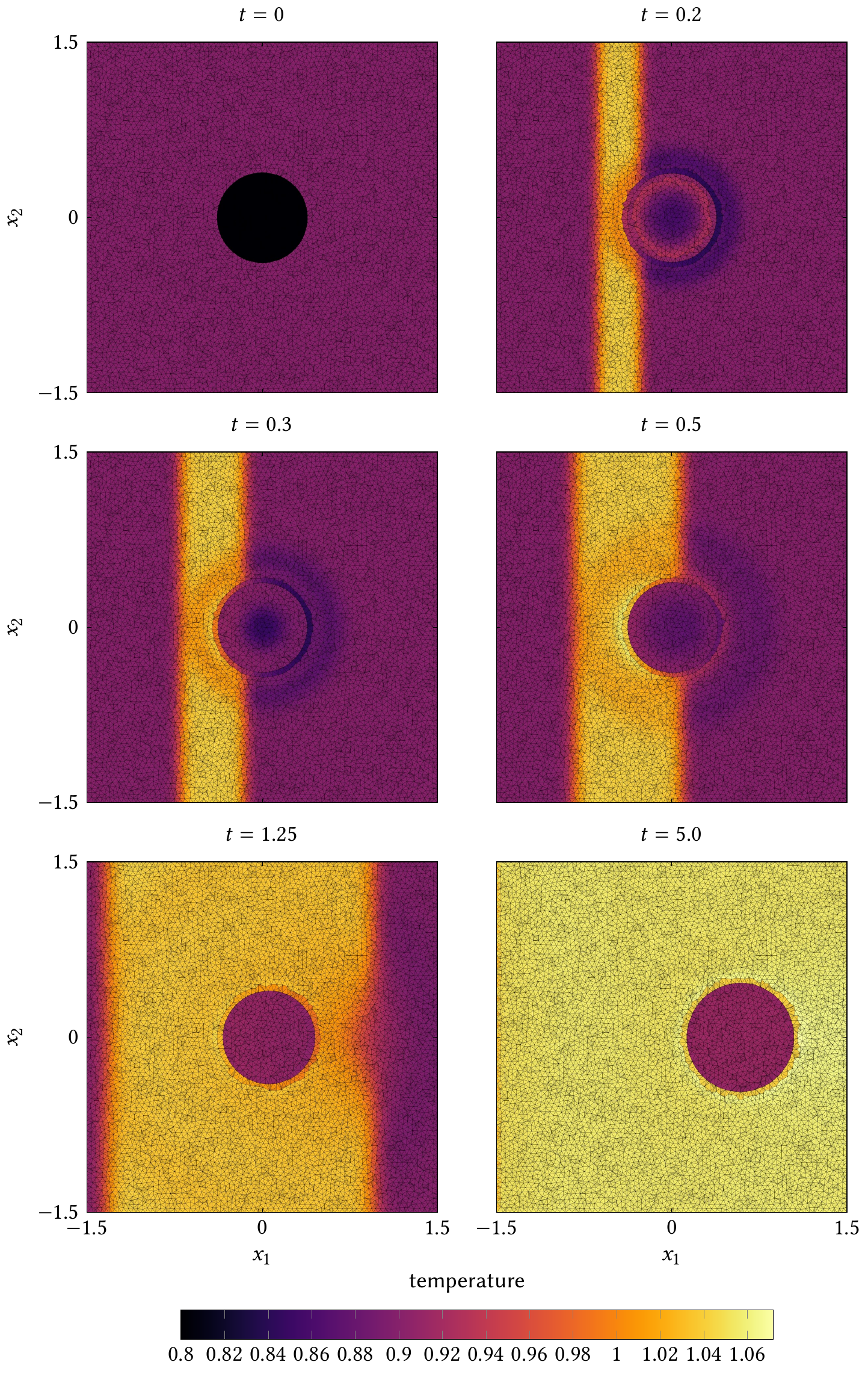}  
\caption{
 Slice through the three-dimensional multiscale simulation of the liquid--vapor flow model at \(x_3 = 0\). 
 The sub-figures depict the temperature of the solution at various time steps. 
 The corresponding plots of the density can be found in Figure~\ref{fig:noniso:shock_droplet:temperature:3d}.
  }
 \label{fig:noniso:shock_droplet:temperature:3d}
\end{figure}

The simulation results are presented in 
Figure~\ref{fig:noniso:shock_droplet:density:3d} for the density and the velocity,
and in Figure~\ref{fig:noniso:shock_droplet:temperature:3d} for the temperature. 
Both figures show two-dimensional slices through the three-dimensional domain at \(x_3 = 0\).
\newline
Initially, the liquid droplet is not in equilibrium with the surrounding vapor, resulting in a compression of the droplet. 
Then the vapor wave hits the droplet and acts upon it, moving the droplet through the computational domain. 
During that phase the droplet slowly grows in volume. 
Compared to the two-dimensional simulation, the deformation is much less pronounced. 
The reason for this is that the fluid has more freedom to move around the droplet, as opposed to the two-dimensional setting.

Finally, let us regard the computational effort of the simulation.\footnote{All computations in this contribution are performed on a desktop computer equipped with an AMD Ryzen Threadripper 2950X 16-core processor and \SI{128}{\giga\byte} RAM.}
Averaged over the 
\num{1e4}
time steps the mesh consists of
\num[scientific-notation=true, round-mode=places, round-precision=1]{2.667542e+06} 
cells with 
\num[scientific-notation=true, round-mode=places, round-precision=1]{3579.837098}
interface facets. 
Each time step takes approximately
\SI[round-mode=places, round-precision=0]{23.8373048}{\second}, 
split into 
\SI[round-mode=places, round-precision=0]{1.810514}{\second} 
for performing the moving mesh operations, and 
\SI[round-mode=places, round-precision=0]{21.58386}{\second}
for computing the finite volume fluxes and updates.
In total (including I/O), this amounts to a runtime of about 
\SI[scientific-notation=false, round-mode=places, round-precision=0]{66.2147356}{\hour}. 
We note however, that the implementation is mostly serial. 
Only the numerical fluxes in the bulk phases are computed in parallel with a shared memory parallelization. 
\newline
Employing the surrogate solver is crucial to run this simulation. 
If we would not have used a surrogate solver, we would need to perform a \mdshort{} simulation at every interface facet. 
This would sum up to 
\SI[scientific-notation=false, round-mode=figures, round-precision=2]{596.639516}{\hour} of (albeit parallelizable) computation time for every time step. 
Compared to that, a single evaluation of the surrogate solver takes only  
\SI[round-mode=places, round-precision=1]{0.156}{\milli\second}. 
Even including the offline phase for generating the surrogate solver (approximately \SI{3000}{\hour}), the additional effort is quickly redeemed after a few time steps.

 \section{Conclusion}
\label{sec:conclusion}

Up to now, there were few methods to solve 
inviscid, temperature-dependent liquid--vapor flow, and the applied methods tend to be either computationally expensive, or are unable to render some physical effects accurately. 
We developed a multiscale model for inviscid liquid--vapor flow, that 
accounts for the interface physics on the molecular scale. 
For this purpose, we utilize the molecular-scale Riemann problem as the link between the continuum scale and the molecular scale.
Thereby, we avoid explicit closures on the continuum scale which are in many cases not valid nor available. 
This enables us to simulate fluid regimes that could not be considered before. 
However, that comes at the cost of expensive \mdshort{} simulations. 
We overcome this bottleneck by employing (mass-conservative) constraint-aware neural networks, rendering the whole implementation efficient enough to perform simulations in two and three space dimensions. 
\newline
Notwithstanding, there remain some open questions.
\begin{itemize}
  \item An important addition to the continuum-scale model and its discretization is to include topological changes such as droplet merging or splitting. 
  For this purpose, sub-models have to be introduced that describe these phenomena in a sharp-interface setting. 
  In the same way, nucleation and vanishing droplets should be supported. 
\item Surface tension effects are an important factor in two-phase flows. 
  However, they are not taken into account in this work.  
  It is straightforward to include them on the continuum scale because the curvature of the phase boundary can be computed directly from the geometric information of the sharp interface in the mesh, see e.g. \cite{popinet:numerical:2018}.
  The local curvature can then be passed to the microscale interface solver. 
  We have already demonstrated this in our work \cite{chalons.magiera.ea:finite:2018}. 
  What is missing is an efficient method to include the curvature effects on the molecular scale. 
\item We consider inviscid flow without heat conduction in the bulk phases.
To include heat and viscous fluxes it suffices to exchange the Euler equations by the Navier--Stokes--Fourier system, while inferring the derivatives of the fluid states from the \mdshort{} simulations. 
\end{itemize}
Regarding the surrogate solver, various improvements can be made. 
Active learning may be utilized and more efficient sampling strategies for generating the training data set can be applied. 
\newline 
The timescale of the \mdshort{} simulations and the size of the interface sampling region are determined heuristically. 
This is due to the lack of rigorous theory for bridging the scales from the molecular scale particle model to a continuum-scale fluid model.

The multiscale model can be naturally extended to the case of more complex fluids, only the \mdshort{} simulations have to be adapted for the corresponding  molecule types. 
Moreover, the multiscale model is easily adaptable to various sharp-interface,  two-phase flow scenarios, while preserving the physical accuracy of \mdshort{} simulations.

\subsection*{Acknowledgements}
\label{sec:acknowledgements}
This work was supported by the Deutsche Forschungsgemeinschaft (DFG, German Research Foundation) through the project  \mbox{SFB--TRR 75} ``Droplet dynamics under extreme ambient conditions'' with the project number 84292822 
and the DFG under Germany's Excellence Strategy - EXC 2075 with the project number 390740016.

\appendix

\section{Elements of Molecular Dynamics Simulations}
\label{appendix:md_elements}

In this section, we explain the tailored  \mdshort{} simulations in detail. 
Starting with the cutoff-potential substituting the full Lennard-Jones potential 
\eqref{eq:md:lennard_jones_potential}.
Subsequently, we give a detailed explanation of other \mdshort{} elements, namely the particular choice of the computational domain, the long-range corrections, and the thermostat.

\paragraph*{Cutoff-Potential}
Solving the equation of motion \eqref{eq:md:equation_of_motion} directly requires the computation of \((\nparticles-1)^2\) pair-interactions.
For a high number of particles this becomes unfeasible. 
To respond to this problem, it is common to apply a localized potential by introducing a cutoff length \(r_{\mathrm{cutoff}} > 0\), above which the potential is set to zero, i.e.
\begin{align} 
  \label{eq:cutoff_potential}
 \phi(r; r_{\mathrm{cutoff}}) \coloneqq \begin{cases}
                         \phi(r), & r < r_{\mathrm{cutoff}}, \\
                         0, & r \geq r_{\mathrm{cutoff}}.
                        \end{cases}
\end{align}
Naturally, using a cutoff-potential comes at a cost: long-range forces are not considered anymore. 
They have to be re-included via long-range corrections or long-range forces, which will be discussed later.

\begin{figure}[htb]
\centering
\subcaptionbox{Standard, cubic computational domain\label{subfig:standard_domain}}[0.45\columnwidth]
{\includegraphics[width=0.175\columnwidth]{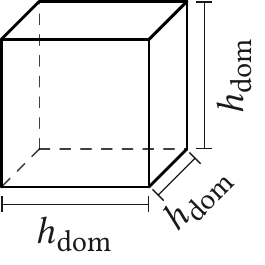}
}
\hspace{0.05\columnwidth}
\subcaptionbox{Elongated cuboid domain with linked-cell approach\label{subfig:cuboid_domain}}[0.45\columnwidth]
{\includegraphics[width=0.4\columnwidth]{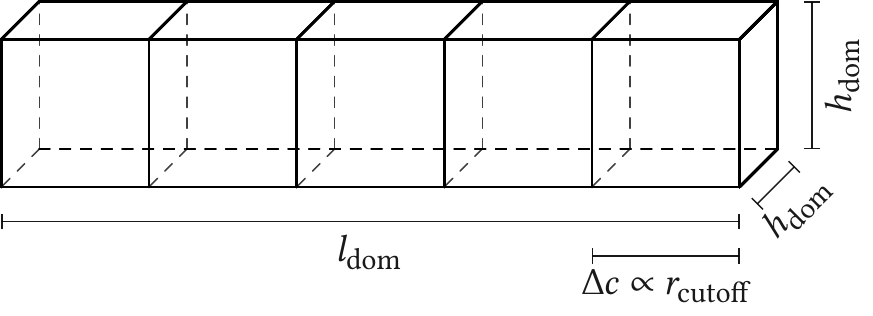}
}
\caption{Two different setups for the computational domain.}
\label{fig:domains}
\end{figure}

\paragraph*{Computational Domain}
\label{sec:md_domain}
We are interested in simulating the fluid dynamics in normal direction of the interface. 
Therefore, it is beneficial to consider an adapted computational domain, that is large in \(x\)-direction and small in the transversal directions,
see Figure~\ref{fig:domains}~\subref{subfig:cuboid_domain}. 
In the transversal directions we consider periodic boundary conditions, whereas in \(x\)-direction reflecting boundary conditions are implemented. 
We fix the height \(h_{\mathrm{dom}}\) in the transversal directions to a multiple of the cutoff-radius \(r_{\mathrm{cutoff}}\) and set the length \(l_{\mathrm{dom}}\) in \(x\)-direction according to the desired fluid density.

\paragraph*{Linked-Cell Approach}
The big advantage of choosing the adapted computational domain is that, in conjunction with the cutoff-potential, we are now able to consider one-dimensional grid cells with cell length \(\Delta c > r_{\mathrm{cutoff}}\) for a more efficient and faster computation of the interaction forces.
This is achieved by computing only the interactions of all particles within a cell and its two neighboring cells. 
In that way, not all pairwise particle distances have to be computed, which saves computational resources. 
This approach is called linked-cell or linked-list approach \cite{allen.tildesley:computer:2017}.
Due to the one-dimensional nature of the grid, we can distribute the particles onto the cells by sorting them according to their \(x\)-coordinate, which is done every  \(\frequencyvar_{\mathrm{sort}} \in \bN\) time steps.

\paragraph*{Long-Range Forces}
Following the approach described in \cite{janecek:long:2006},  
we compute the forces depending on the cell-wise local density average. 
Therefore, the cell-averaged densities \(\overline{\rho}_{C_k}\) are computed at each point in time for every cell \(C_k\). 
Then, a particle \(i\) in cell \(C_k\) is subject to the long-range force 
\begin{align} 
  \label{eq:longrange_force}
  \vf_{C_k, \mathrm{LR}} = \sum_{l \neq k-1, k, k+1} 
  2 \pi (c_k - c_l)  \phi(\abs{c_k - c_l}) \overline{\rho}_{C_k} \Delta c
  \quad \text{for } \vx_i \in C_k,
\end{align}
where \(c_l\) denotes the center point of cell \(C_l\).
Note that the long-range forces do not depend on the particle position and can therefore be precomputed for each cell.

\paragraph*{Statistical Averaging}
\label{sec:md_averaged_quantities}
By solving the \mdshort{} system \eqref{eq:md:equation_of_motion}, microscale information about the fluid system is generated. 
In most situations, the trajectory of every single particle is not of interest. 
Instead, the focus lies on the continuum-scale behavior of the fluid.
To convert the microscale dynamics to continuum-scale quantities, we apply methods from the domain of statistical mechanics. 
By employing the Irving--Kirkwood formulas \cite{irving.kirkwood:statistical:1950} it is possible to formulate microscopic distributions of some key quantities.
Thus, the microscopic instantaneous density \(\widehat{\rho}(x,\mdtimevar)\) and momentum \(\widehat{\vm}(x,\mdtimevar)\) distributions are given by
\begin{align} \label{eq:irving-kirkwood:md}
\begin{aligned}
 \widehat{\rho}(x,\mdtimevar) & = \sum_{i=1}^{\nparticles} m_i \, \delta_{x_i(\mdtimevar)}(x), 
 &
 \widehat{\vm}(x,\mdtimevar) & = \sum_{i=1}^{\nparticles} m_i \, \vv_i(\mdtimevar) \, \delta_{x_i(\mdtimevar)}(x),     
 \end{aligned}
\end{align}
where \(m_i > 0\), \(x_i \in \bR\), \(\vv_i \in \bR^3\) are the mass, position in \(x\)-direction and velocity of the \(i\)-th particle.
Here, \(\delta_{x_i(\mdtimevar)}(\,\cdot\,)\) denotes the Dirac delta distribution centered at \({x_i}\).
The local temperature distribution can be formally understood as
\begin{align}
  \label{eq:md_temperature_distribution}
 \widehat{T}(x,\mdtimevar) \coloneqq \frac{1}{d} \left( \sum_{i=1}^{\nparticles} \delta_{x_i(\mdtimevar)}(x) \right)^{-1} \sum^{\nparticles}_{i = 1} (m_i (\vv_i - \overline{\vv}) \cdot (\vv_i - \overline{\vv})) \, \delta_{x_i(\mdtimevar)}(x).
\end{align}
However, the temperature is only defined as an averaged quantity and is undefined in the absence of particles.

\paragraph*{Thermostat}
To provide meaningful \mdshort{} simulations, it is important to control the temperature. 
Typically, this is done via a so-called thermostat during the simulations.
We will apply one of the most basic thermostats: the velocity rescaling thermostat, see e.g. \cite{bussi.donadio.ea:canonical:2007}. 
It controls the temperature by rescaling the 
relative velocity 
\((\vv_i - \overline{\vv})\) 
of every particle by
\begin{equation}
 \lambda = \sqrt{T_{\mathrm{ref}} / T_{\mathrm{inst}}},
\end{equation}
where \(T_{\mathrm{ref}} \geq 0\) is the prescribed target temperature and \(T_{\mathrm{inst}}\) is the measured, instantaneous temperature, averaged over the whole system using \eqref{eq:md_temperature_distribution}.

\paragraph*{\mdshort{} Time-Step Algorithm}

The following algorithm describes the computation of a single time step of the \mdshort{} simulation.
We use the Velocity Verlet algorithm \cite{verlet:computer:1967} to solve \eqref{eq:md:equation_of_motion} for a time step from \(\mdtimevar_n = n \, \Delta \mdtimevar\) to \(\mdtimevar_{n+1}\), with \(\Delta \mdtimevar > 0\) and \(n \in \bN\).

\begin{boxalgorithm}[label=alg:md_step]{Molecular Simulation Time-Step with Linked Cells}
\textbf{Input:}
particle state arrays \(\vx[i] = \vx_i(\mdtimevar_n)\), \(\vv[i] = \vv_i(\mdtimevar_n)\), and \(\va[i] = \va_i(\mdtimevar_n)\) for \(i = 1, \ldots, \nparticles\) at \(\mdtimevar_n \geq 0\). 

\textbf{Parameters:}
time step \(\Delta \mdtimevar > 0\),
  sorting frequency \(\frequencyvar_{\mathrm{sort}} \in \bN\),
  linked-cell size \(\Delta c > 0\).
\tcbsubtitle{Algorithm}
\begin{itemize}
  \item If \(n \mod \frequencyvar_{\mathrm{sort}} = 0\): sort particles w.r.t. \(x\)-coordinates and store particles indices for each cell \(C_k\). 
\item Compute cell-averaged densities \(\overline{\rho}_{C_k}\) and the long-range forces \(\vf_{C_k, \mathrm{LR}}\) for each cell \(C_k\), see \eqref{eq:longrange_force}.  
  \item Update positions, for all \(i = 1, \ldots, \nparticles\),
  \begin{align}    
    \vx[i] &\pluseq \Delta \mdtimevar \, \vv[i]  + \tfrac{1}{2} \Delta \mdtimevar^2 \, \va[i].  
  \end{align}  
  \item Apply boundary conditions.
\item Compute intermediate velocities, for all \(i = 1, \ldots, \nparticles\)
  \begin{align}
      \vv[i] &\pluseq \tfrac{1}{2} \Delta \mdtimevar \, \va[i].
  \end{align}
  \item Update accelerations, for all \(i = 1, \ldots, \nparticles\),
  \begin{align}
    \va[i] = \tfrac{1}{m_i} (\vf_i + \vf_{C_k, \mathrm{LR}}),
  \end{align} 
 using the pair-interaction potential \(\Phi\), the (updated) positions \(\vx[\cdot]\), and equations \eqref{eq:cutoff_potential}, 
\eqref{eq:longrange_force}.
  \item Update velocities, for all \(i = 1, \ldots, \nparticles\),
  \begin{align}
  \vv[i] \pluseq  \tfrac{1}{2} \Delta \mdtimevar \, \va[i].
  \end{align}
\end{itemize}
\textbf{Result:} \(\vx_i(\mdtimevar_{n+1}) = \vx[i]\), \(\vv_i(\mdtimevar_{n+1}) = \vv[i]\), and \(\va_i(\mdtimevar_{n+1}) = \va[i]\).
\end{boxalgorithm}
 
\paragraph*{Molecular-Scale Riemann Problem}
The following algorithm describes the initialization of the molecular-scale Riemann problem in detail. 

\begin{boxalgorithm}[label=alg:md_riemann_init]{Molecular-Scale Riemann Problem Initialization}  
  \textbf{Input:}
  continuum fluid states \((\rho_\phasel, v_\phasel, T_\phasel)\), \((\rho_\phaser, v_\phaser, T_\phaser)\). 
  
  \textbf{Parameters:}
  total number of particles \(\nparticles\),
  number of thermalization steps \(N_{\mathrm{therm. steps}}\),
  thermostat frequency \(\frequencyvar_{\mathrm{th}}\),
  domain width factor \(w\),
  maximum particle ratio \(\alpha_{\mathrm{max}}\),
  domain length ratio \(\alpha_{\mathrm{dom}}\).
\tcbsubtitle{Algorithm}
\begin{itemize}
  \item Particle ratio: \(\alpha =(1 + \alpha_{\mathrm{dom}} \rho_\phaser / \rho_\phasel)^{-1}\), clamped at \(1-\alpha_{\mathrm{max}}\) and \(\alpha_{\mathrm{max}}\).
  \item Number of particles in \(\phaserl\)-phase: \(N_\phaserl = \alpha \nparticles\).
  \item Domain size in transversal directions: \(l_{yz} = w \sigma r_{\mathrm{cutoff}}\).
  \item Volume of each phase-subdomain \(\mddomain_\phaserl\): \(\mathrm{vol}_\phaserl = \frac{m_0 N_\phaserl}{\rho_\phaserl}\).
  \item Domain width for each phase subdomain \(\mddomain_\phaserl\) in \(x\)-direction: \(l_\phaserl = \frac{\mathrm{vol}_\phaserl}{l_{yz}}\). 
  \item Interface gap to avoid overlapping particles: \(l_{\mathrm{gap}} = 2^{\frac{1}{6}} \sigma\).
  \item Define the domains
  \begin{align}
    \mddomain_\phasel &\coloneqq [0, l_\phasel] \times [0, l_{yz}]^2, 
    &
    \mddomain_{\phasel, \mathrm{gap}} &\coloneqq [0, l_\phasel - \tfrac{1}{2} l_{\mathrm{gap}}] \times [0, l_{yz}]^2, \\
    \mddomain_\phaser &\coloneqq [l_\phasel, l_\phaser + l_\phasel] \times [0, l_{yz}]^2,  
    &
    \mddomain_{\phaser, \mathrm{gap}} &\coloneqq [l_\phasel + \tfrac{1}{2} l_{\mathrm{gap}}, l_\phaser + l_\phasel] \times [0, l_{yz}]^2.
  \end{align} 
  \item Distribute \(N_\phaserl\) particles on a uniform lattice in \(\mddomain_{\phaserl, \mathrm{gap}}\). 
  If particles remain, distribute them randomly, while avoiding overlaps (using Algorithm~\ref{alg:poisson_disc}). 
  \item \emph{Thermalization phase:}
  run \(N_{\mathrm{therm. steps}}\) of a \mdshort{} simulation in each subdomain \(\mddomain_{\phaserl}\), and each \(\frequencyvar_{\mathrm{th}}\)-th step apply a thermostat with target temperature \(T_{\phaserl}\). \item Remove average velocities and add continuum fluid velocities \(v_\phaserl\) to all particles in \(\mddomain_{\phaserl}\). 
  \item Combine domains \(\mddomain \coloneqq \mddomain_{\phasel} \cup \mddomain_{\phaser}\).
\end{itemize}
\textbf{Result:} \(\nparticles\) particles in a microscale Riemann problem configuration with corresponding continuum fluid states \((\rho_\phasel, v_\phasel, T_\phasel)\) and \((\rho_\phaser, v_\phaser, T_\phaser)\). 
\end{boxalgorithm}

\paragraph*{Molecular-Scale Interface Tracking}
The following algorithm describes the tracking of the phase boundary on the molecular scale.

\begin{boxalgorithm}[label=alg:md_tracking]{Molecular-Scale Interface Tracking}
\textbf{Input:} 
particle positions \(\vx_i(\mdtimevar_n)\), \(i = 1, \ldots, \nparticles\), at time \(\mdtimevar_n >0\), \(n \in \bN\);
former interface positions \(\mdinterface(t_{j})\) for \(j = \max(n - n_{\mathrm{buffer}}, 0), \ldots, \max(n-1, 0)\).

\textbf{Parameters:}  
RBF kernel parameter \(\gamma_{\mathrm{RBF}}\),
interface tracking width ratio \(\alpha_{\mathrm{track}}\),
buffer size \(n_{\mathrm{buffer}}\).

\tcbsubtitle{Algorithm}
\begin{itemize}
  \item Get the new interface position \(\mdinterface(t_{n})\) by  
  solving the optimization problem in 
  \(B_{w_{\mathrm{track}}}(\mdinterface(\mdtimevar_{n-1})) \allowbreak \coloneqq  \allowbreak (\mdinterface(\mdtimevar_{n-1}) - w_{\mathrm{track}}, \mdinterface(\mdtimevar_{n-1}) + w_{\mathrm{track}}) \allowbreak \subset \bR\), 
  with \(w_{\mathrm{track}} \coloneqq \alpha_{\mathrm{track}} l_{\mathrm{dom}}\),
\begin{align}
  \mdinterface(\mdtimevar_n) = \argmax_{x \in B_{l_{\mathrm{offset}}}(\mdinterface(\mdtimevar_{n-1})) }  \abs{\partial_x \widehat{\rho}(x,\mdtimevar_n)}.
\end{align}
\item If \(n \geq 2\), perform linear regression over the data \((\mdtimevar_j, \mdinterface(\mdtimevar_j))\), \(j = \max(n - n_{\mathrm{buffer}}, 0), \ldots, n\).
The slope \(s\) of the fitted line is the averaged interface speed \(s(\mdtimevar_n) = s\).
\end{itemize}
\textbf{Result:}   
interface position \(\mdinterface(\mdtimevar_n)\) and interface speed \(s(\mdtimevar_n)\).
\end{boxalgorithm}

\section{Sampling Algorithms}
\label{appendix:sampling_algorithms}

\paragraph*{Uniform Sampling in a Convex Domain}
We seek to draw random samples inside a convex domain \(\Omega \subset \bR^d\).
Commonly, this is done via rejection sampling, but this can become very expensive in high-dimensional spaces. 
Alternatively, we perform uniform sampling on a mesh of the convex domain. 
Here we exploit the fact that the Dirichlet-distribution yields a uniform point distribution after projection on a simplex. 

\begin{boxalgorithm}[label=alg:convex_hull_sampling]{Convex Set Sampling}
\textbf{Input:}  
point cloud \(P = \{\vp_0, \ldots, \vp_M\} \subset \bR^d\), 
number of new samples \(N\).
\tcbsubtitle{Algorithm}
\begin{itemize} 
  \item Construct convex hull of \(P\).
  \item Construct \(d\)-dimensional Delaunay mesh of the points on the convex hull.
  \item Calculate the volume of each simplex in the mesh.
  \item For \(i = 1, \ldots, N\) do
\begin{itemize}
\item \(T_j \leftarrow\) choose a random simplex from the Delaunay mesh with probability proportional to the simplex volume.
  \item \(\vec{\alpha} \leftarrow\) Get a \((d+1)\)-dimensional vector drawn from a uniform Dirichlet distribution.
  \item \(\vx_{i} \leftarrow\) linear combination of the vertices of \(T_j\) with the coefficients of \(\vec{\alpha}\).
  \end{itemize}
\end{itemize}
\textbf{Result:} \(\{\vx_1, \ldots, \vx_N\}\).
\end{boxalgorithm}

\paragraph*{Distance-Maximizing Sampling}
Uniform random sampling might result in bad space-filling behavior, especially for sparse data in high-dimensional spaces. 
We use distance-maximizing sampling using a multidimensional variant of the best-candidate algorithm of \cite{mitchell:spectrally:1991}. 
However, by using the standard euclidean metric to measure the distance this can lead to a high concentration of samples near the domain boundary \cite{santner.williams.ea:design:2003}. 
Instead, we employ the Minkowski distance 
\begin{align}
  \operatorname{dist}(\vx,\vy) \coloneqq \left( \sum^{d}_{i=1} \abs{x_i - y_i}^p \right)^{\frac{1}{p}}
\end{align}
with \(p \in (0,1)\). 
This choice appears to provide a much better performance in high dimensional spaces \cite{aggarwal.hinneburg.ea:surprising:2001} and avoids a high concentration of the points near the domain boundary.
As a rule of thumb, we set \(p = 4^{-d}\). 
This empirical choice for \(p\) appears to produce the desired results.

\begin{boxalgorithm}[label=alg:poisson_disc]{Distance-Maximizing Sampling}
\textbf{Input:}  
convex sampling domain \(\Omega \subset \bR^d\), 
number of new samples \(N \in \bN\), 
number of test samples \(K \in \bN\).
\tcbsubtitle{Algorithm}
\begin{itemize} 
  \item Sample one uniform random point in \(\Omega\) using Algorithm~\ref{alg:convex_hull_sampling} to initialize \(X = \{\vx_1\}\).
  \item For \(i = 2, \ldots, N\) do 
  \begin{itemize}
\item \(Y = \{\vy_1, \ldots, \vy_K\} \leftarrow\) sample \(K\) uniform random points in \(\Omega\) using Algorithm~\ref{alg:convex_hull_sampling}.
    \item \(j \leftarrow \operatorname{arg\,max}_{k = 1, \ldots,K} \operatorname{dist}(\vy_k, X)\).
    \item \(\vx_i \leftarrow \vy_j\).
    \item \(X \leftarrow X \cup \{\vx_i\}\).
  \end{itemize}
\end{itemize}
\textbf{Result:} 
Samples \(X = \{\vx_1, \ldots, \vx_N\}\).
\end{boxalgorithm}

 \section{Surrogate Solver Details for the Liquid--Vapor Flow Model}
\label{appendix:surrogate_details_euler}

In this section, we provide some details regarding the surrogate solver for liquid--vapor flow discussed in Section~\ref{sec:multiscale_nonisothermal}, 
namely the data set generation and the outline of the training procedure.

\paragraph*{Data Set Generation}
\label{sec:data_generation:euler}
To be able to train a surrogate solver,  an appropriate data set has to be generated. 
The general outline of the procedure is described in Section~\ref{sec:surrogate}.  
First, we have to define the model-specific data range for the input data set \(\datasetvar_{\mathrm{in}}\). 
Each input data point is given by a tuple \((\rho_\phasel, v_\phasel, T_\phasel,  \rho_\phaser,  v_\phaser,  T_\phaser)\), with \(v_\phasel \equiv 0\) after transforming the reference frame.  
The range of the \(\phaser\)-phase densities \(\rho_\phaser\) and temperatures \(T_\phaser\) is defined by the convex set formed by the following points
\begin{align}
  (\rho_{\mathrm{min},\phaser}, T_{\mathrm{min}}),
  (\rho_{\mathrm{min},\phaser}, T_{\mathrm{max}}),
  (\rho_{\mathrm{max},\phaser}, T_{\mathrm{min}}),
  (\rho_{\mathrm{c}}, T_{\mathrm{c}}),
  (\rho_{\mathrm{c}}, T_{\mathrm{max}}), 
\end{align}
with 
\(\rho_{\mathrm{min},\phaser} = 10^{-4}\), 
\(\rho_{\mathrm{max},\phaser} = \num{0.02}\), 
\(\rho_{\mathrm{c}} = \num{0.31}\), 
\(T_{\mathrm{min}} = \num{0.4}\), 
\(T_{\mathrm{c}} = \num{1.32}\), 
and 
\(T_{\mathrm{max}} = \num{2.5}\). 
Analogously, the \(\phasel\)-phase range for the densities \(\rho_\phasel\) and temperatures \(T_\phasel\) is given by the convex set formed by the points
\begin{align}
  (\rho_{\mathrm{min},\phasel}, T_{\mathrm{min}}),
  (\rho_{\mathrm{c}}, T_{\mathrm{c}}),
  (\rho_{\mathrm{max},\phasel}, T_{\mathrm{min}}),
  (\rho_{\mathrm{max},\phasel}, T_{\mathrm{c}}),
  (\rho_{\mathrm{c}}, T_{\mathrm{max}}),
\end{align}
with 
\(\rho_{\mathrm{min},\phasel} = \num{0.7}\), and  
\(\rho_{\mathrm{max},\phasel} = \num{1.0}\).
Finally, the \(\phaser\)-velocity lies in the interval \(v_\phaser \in (v_{\mathrm{min}}, v_{\mathrm{max}})\), with 
\(v_{\mathrm{min}} = -2.5\), 
and
\(v_{\mathrm{max}} = 1.5\).
Altogether, this gives the definition of the input bounding domain \(B_{\mathrm{in}}\).
\newline
Inside the input bounding domain, \(N_{\mathrm{in}} = \num{6000}\) data points are sampled, using Algorithm~\ref{alg:poisson_disc}, while exploiting that \(v_{\phasel} \equiv 0\), i.e. the sampling domain is five-dimensional.
This results in the initial data set 
\(\datasetvar_{\mathrm{in}} = \{ (\rho_\phasel, \allowbreak v_\phasel \equiv 0, \allowbreak T_\phasel, \allowbreak \rho_\phaser, \allowbreak v_\phaser, \allowbreak T_\phaser)_i 
\setsep  
i = 1,\ldots,N_{\mathrm{in}}\}\).
\newline
For each \((\rho_\phasel, v_\phasel \equiv 0, T_\phasel, \rho_\phaser,  v_\phaser,  T_\phaser) \in \datasetvar_{\mathrm{in}}\) we evaluate the  \mdshort{} interface solver (Algorithm~\ref{alg:md_riemann}), i.e.
\begin{align}
  \label{eq:md:data_evaluation}
  \rmd(\rho_\phasel, v_\phasel \equiv 0, T_\phasel, \rho_\phaser,  v_\phaser,  T_\phaser) 
  &= 
  (\rho^*_\phasel, m^*_\phasel, T^*_\phasel, \rho^*_\phaser, m^*_\phaser, T^*_\phaser, s). 
\end{align}
We gather the valid output in the output data set 
\(\datasetvar_{\mathrm{out}} 
= \{ (\rho^*_\phasel, \allowbreak m^*_\phasel, \allowbreak T^*_\phasel, \allowbreak \rho^*_\phaser, \allowbreak m^*_\phaser, \allowbreak T^*_\phaser, \allowbreak s)_i 
\setsep 
i = 1, \allowbreak \ldots, \allowbreak N_{\mathrm{data}}\}\)
 and identify each output sample with its corresponding input 
\((\rho_\phasel\), \(v_\phasel \equiv 0\), \(T_\phasel\), \(\rho_\phaser\), \(v_\phaser\),  \(T_\phaser)\), 
resulting in the complete data set
\begin{align}
  \datasetvar = \Bigl\{ 
  \bigl((\rho_\phasel, v_\phasel \equiv 0, T_\phasel, \rho_\phaser,  v_\phaser,  T_\phaser), 
 (\rho^*_\phasel, m^*_\phasel, T^*_\phasel, \rho^*_\phaser, m^*_\phaser, T^*_\phaser, s)\bigr)_i \setsep  i = 1,\ldots,N_{\mathrm{data}}
 \Bigr\}.
\end{align}

\begin{remark}[Repeating Molecular Dynamics Simulations]
  \label{remark:repeating:iso_md}
  As the \mdshort{} simulations are non-deterministic and subject to noise, we repeat the \mdshort{} simulation, i.e. the evaluation of \eqref{eq:md:data_evaluation}, for each data point \num{3} times. 
  This results in a larger data set and gives us a rough indicator of the accuracy of the \mdshort{} simulations. Which --- implicitly --- can help us to avoid overfitting during the neural network training.
\end{remark}

\begin{figure}[tp]
\centering
\includegraphics[width=0.95\columnwidth]{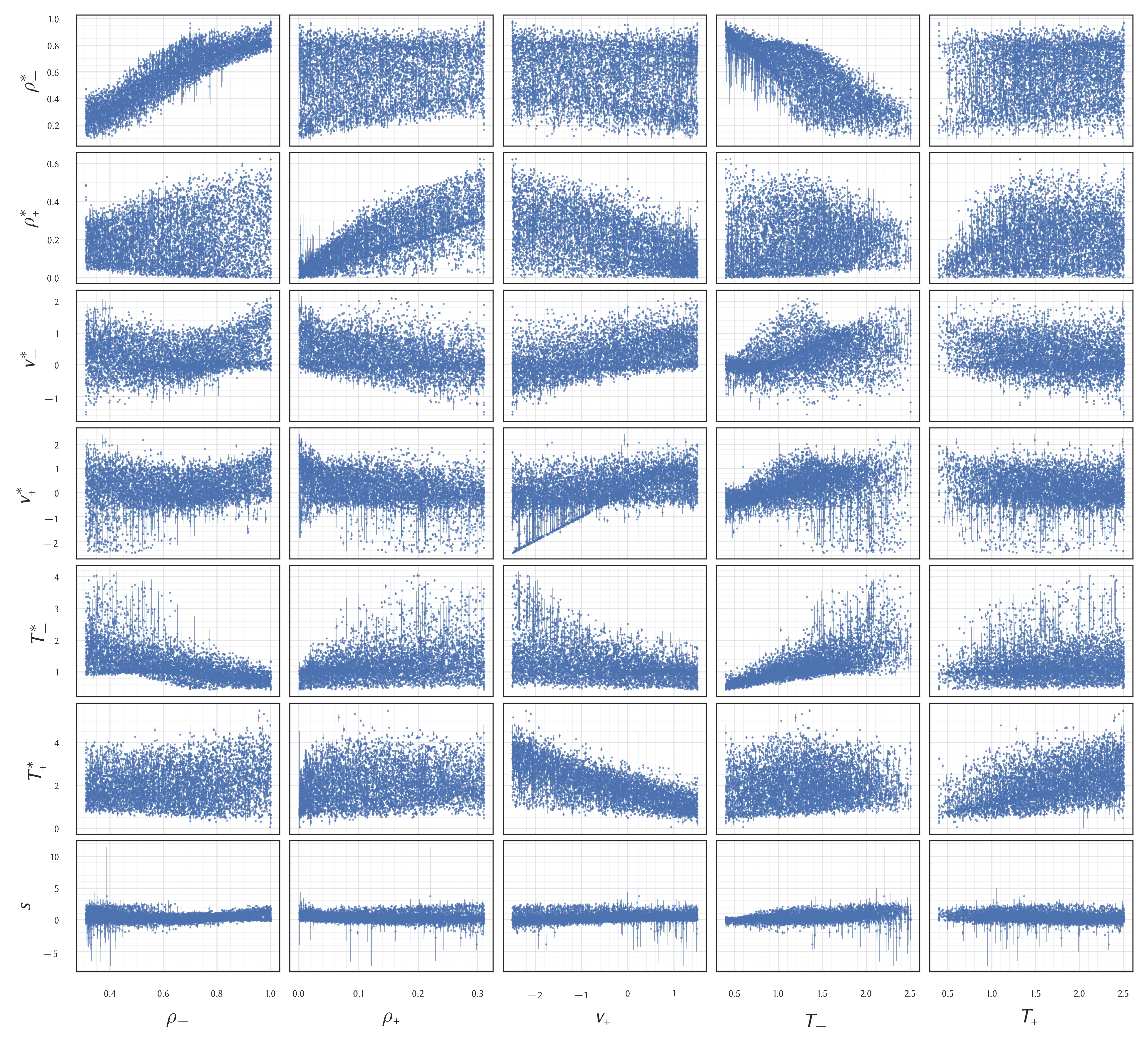}
\caption{The data set for the liquid--vapor flow multiscale model. 
The data set consists of \num{17865} data points (with \num{5955} unique initial states).
The scatter plots show the relation of the input states \(\rho_\phasel\), \(T_\phasel\), \(\rho_\phaser\), \(v_\phaser\), \(T_\phaser\) to the wave states \(\rho^*_\phasel\), \(v^*_\phasel\), \(T^*_\phasel\), \(\rho^*_\phaser\), \(v^*_\phaser\), \(T^*_\phaser\), and wave speed \(s\). 
The error-bars show the data range of each separate \mdshort{} simulation.
The input variable \(v_{\phasel}\) is not shown as it equals zero throughout the data set.  
}
\label{fig:noniso_dataset:errorbars}
\end{figure}

The data set is available at \cite{magiera:data:2021}.
To get an overview of the data set, the resulting data set is illustrated in Figure~\ref{fig:noniso_dataset:errorbars}, 
where we plot the output values over the input values, including the range of the output values, which is indicated by error bars.

A single run of a single \mdshort{} simulation (i.e.  Algorithm~\ref{alg:md_riemann}) takes between \num{9} and \SI{11}{\minute}. 
For the whole data set this amounts to roughly \SI{3000}{\hour} of computing time. 
The sampling of the data points is however fully parallelizable, and therefore, the workload can be easily split among numerous machines.

\paragraph*{Neural Network Training}
\label{sec:network_training:euler}

\begin{figure}[tp]
\centering
\includegraphics[width=0.45\columnwidth]{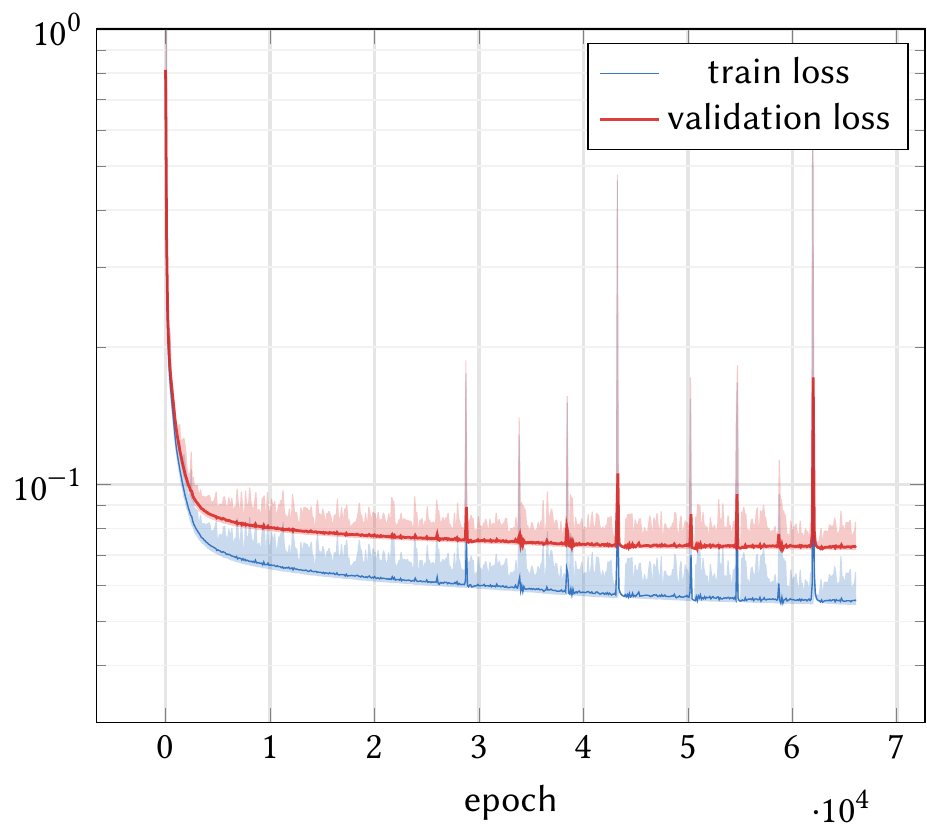}
\caption{Evolution of the loss on the training and the validation data set during the training of the neural network surrogate solver for the liquid--vapor flow multiscale model.}
\label{fig:noniso_nn_training_history}
\end{figure}

The training of the neural network proceeds in the usual way and follows our approach used in \cite{magiera.ray.ea:constraint:2020} for \cres{}-networks. 
For training and employing  we use the neural network framework PyTorch \cite{paszke.gross.ea:pytorch::2019}.
The particular parameters for this model are outlined in the following. 
We use the data set \(\datasetvar\) of \ref{sec:data_generation:euler} and split it into  a training data set \(\datasetvar_{\mathrm{train}}\) with \num{14889} samples and a validation data set \(\datasetvar_{\mathrm{val}}\) with \num{2976} samples.
We ensure that all repetitions of the \mdshort{} simulations (see Remark~\ref{remark:repeating:iso_md}) belong to either the training or the validation data set. 
The neural network consists of \num{5} hidden layers with \num{50} nodes each. 
The learning rate of the Adam optimizer \cite{kingma.ba:adam:2014} is set to \num{e-3}, the weight decay parameter \(\alpha_{\mathrm{wd}}\) to \num{e-3}, and the gradient norm clipping parameter to \num{e3}.  
The training and validation errors during the training procedure are shown in Figure~\ref{fig:noniso_nn_training_history}.

\section{Temperature-Dependent Multiscale Model: Parameter Tables}
\label{appendix:parameter_table:noniso}

\begin{table}[H]
  \centering
\begin{tabular}{lr}
    \toprule
\multicolumn{2}{c}{\textbf{Molecular Dynamics}} \\ \midrule
    \multicolumn{2}{l}{\textit{Lennard--Jones potential}} \\
    $\sigma$ & \num{1} \\
    $\varepsilon$ & \num{1} \\
    cutoff-radius 
    $r_{\mathrm{cutoff}}$ & \num{2.5} \\
    particle mass $m_0$ & \num{1} \\ 
    \midrule
    \multicolumn{2}{l}{\textit{Simulation}} \\
    number of particles
    $\nparticles$ &  \num{32768} \\
number of time steps
    $N_{\mathrm{end}}$ & \num{e3} \\
    processing frequency
    $\frequencyvar_{\mathrm{pr}}$ & \num{100} \\    
    time step
    $\Delta t$ & \num{e-3}  \\
domain width factor 
    $w$ & \num{3.5} \\  
    domain length ratio 
    $\alpha_{\mathrm{dom}}$ & \num{3.0} \\ 
    maximum particle ratio 
    $\alpha_{\mathrm{max}}$ & \num{0.99}  \\
linked-cell size 
    $\Delta c$ & \num{1.2} $\sigma r_{\mathrm{cutoff}}$ \\
    sorting frequency 
    $\frequencyvar_{\mathrm{sort}}$ & \num{20} \\
thermalization steps
    $N_{\mathrm{therm. steps}}$ & \num{5e2} \\
    thermalization frequency
    $\frequencyvar_{\mathrm{th}}$  & \num{100} \\      
particles per thermostat cell $n_{\mathrm{therm,L}}$ & \num{500} \\ 
    particles per thermostat cell $n_{\mathrm{therm,R}}$ & \num{50} \\   
time sampling ratio 
    $\alpha_{\mathrm{t-smpl}}$ & \num{0.2} \\
    sampling width 
    $w_{\mathrm{sr}}$ & \num{50} \\
    sampling interface offset 
    $o_{\mathrm{sr}}$ & \num{2.5} \\
    interface tracking buffer size 
    $n_{\mathrm{buffer}}$ & \num{5} \\
    interface tracking width 
    $\alpha_{\mathrm{track}}$ & \num{0.1} \\
    RBF parameter
    $\gamma_{\mathrm{RBF}}$ & \num{5e-4} \\ 
\bottomrule
  \end{tabular}
\caption{Parameter table for the molecular dynamics simulation in case of the temperature-dependent, two-phase flow multiscale model.}
  \label{tab:noniso_paramater}  
\end{table}

\begin{table}[H]
  \centering
\begin{tabular}{lr}
    \toprule
\multicolumn{2}{c}{\textbf{\ipfv{} 1d}} \\ \midrule
end time $t_{\mathrm{end}}$ & \num{0.1} \\
    time step $\Delta t$ & \num{e-5} \\
domain $\Omega$ & $(-1, 1)$ \\
    base cell width $\Delta x$ & \num{e-3} \\
    $\Delta x_{\mathrm{min}}$ & $\sfrac{1}{2} \Delta x$ \\    
    $\Delta x_{\mathrm{max}}$ & $\sfrac{3}{2} \Delta x$ \\  
Lax--Friedrichs $\alpha_{\mathrm{LF}}$ & \num{1.0} \\
\bottomrule
  \end{tabular}
\caption{Parameter table for the one-dimensional \ipfv{}-method in case of the liquid--vapor flow multiscale model.}
  \label{tab:noniso_paramater:mmfv:1d}  
\end{table}

\begin{table}[H]
  \centering
\begin{tabular}{lr}
    \toprule
\multicolumn{2}{c}{\textbf{\ipfv{} 2d}} \\ \midrule
time step $\Delta t$ & \num{1e-4} \\
domain $\Omega$ & $(-1.5, 1.5)^2$ \\
    max edge length $\Delta x$ & \num{1e-2} \\
motion regularization $\lambda_{\mathrm{motion}}$ & \num{1e-4} \\    
Lax--Friedrichs $\alpha_{\mathrm{LF}}$ & \num{0.5} \\
\bottomrule
  \end{tabular}
\caption{Parameter table for the two-dimensional \ipfv{}-method in case of the liquid--vapor flow multiscale model.}
  \label{tab:noniso_paramater:mmfv:2d}  
\end{table}

\begin{table}[H]
  \centering
\begin{tabular}{lr}
    \toprule
\multicolumn{2}{c}{\textbf{\ipfv{} 3d}} \\ \midrule
time step $\Delta t$ & \num{5e-4} \\
domain $\Omega$ & $(-1.5, 1.5)^3$ \\
    max edge length $\Delta x$ & \num{3.5e-2} \\
motion regularization $\lambda_{\mathrm{motion}}$ & \num{1e-4} \\  
    sliver threshold $\vartheta_{\mathrm{sliver}}$  & \SI{12}{\degree} \\ 
Lax--Friedrichs $\alpha_{\mathrm{LF}}$ & \num{0.5} \\
\bottomrule
  \end{tabular}
\caption{Parameter table for the three-dimensional \ipfv{}-method in case of the liquid--vapor flow multiscale model.}
  \label{tab:noniso_paramater:mmfv:3d}  
\end{table}

\printbibliography

\end{document}